\documentclass[a4, amsmath, amsthm, amssymb]{amsart}

\usepackage{xcolor}
\usepackage{graphicx}
\numberwithin{equation}{section}

\newtheorem{Lemma}{Lemma}[section]
\newtheorem{Theorem}[Lemma]{Theorem}
\newtheorem{Proposition}[Lemma]{Proposition}
\newtheorem{Corollary}[Lemma]{Corollary}
\newtheorem{Assumption}[Lemma]{Assumption}

\newtheorem{remark}[Lemma]{Remark}
\newtheorem{definition}[Lemma]{Definition}
\newtheorem{example}[Lemma]{Example}
\newtheorem{exercise}[Lemma]{Exercise}

\def\bt{\begin{Theorem}}
\def\et{\end{Theorem}}
\def\bl{\begin{Lemma}}
\def\el{\end{Lemma}}
\def\bp{\begin{Proposition}}
\def\ep{\end{Proposition}}
\def\bcor{\begin{Corollary}}
\def\ecor{\end{Corollary}}
\def\brem{\begin{remark}\rm }
\def\erem{\hfill $\diamondsuit$\end{remark}}

\def\bpf{\begin{proof}}
\def\epf{\end{proof}}

\def\bedef{\begin{definition}\rm }
\def\endef{\hfill $\diamondsuit$\end{definition}}
\def\beg{\begin{example}\rm}
\def\eeg{\hfill $\diamondsuit$\end{example}}
\def\bex{\begin{exercise}\rm}
\def\eex{\hfill $\diamondsuit$\end{exercise}}

\def\berem{\begin{remark}\rm}
\def\erem{\hfill $\diamondsuit$\end{remark}}

\def\bq{\begin{quote}}
\def\eq{\end{quote}}

\def\bc{\begin{center}}
\def\ec{\end{center}}
\def\noi{\noindent}
\def\vsh{\vskip .5cm}

\def\vsq{\vskip .25cm}
\def\beq{\begin{equation}}
\def\eeq{\end{equation}}
\def\beqarray{\begin{eqnarray*}}
\def\eeqarray{\end{eqnarray*}}
\def\<{\leftangle}
\def\>{\rightangle}
\def\({\left(}
\def\){\right)}

\def\<{\langle}
\def\>{\rangle}
\def\q{\quad}

\def\a{\alpha}

\def\g{\gamma}

\def\d{\delta}

\def\t{\tau}

\def\e{\varepsilon}

\def\z{\zeta}

\def\h{\hbox}

\def\w.r.t.{with respect to}
\def\R{{\mathbb{R}}}

\def\G{\Gamma}

\def\N{{\mathbb{N}}}

\def\O{\Omega}
\def\o{\omega}

\def\ben{\begin{enumerate}}
\def\een{\end{enumerate}}
\def\i{\item}
\def\bit{\begin{itemize}}
\def\eit{\end{itemize}}

\begin{document}

\title[Nonlinear
inverse problem]{A new regularization method for a parameter identification problem in a non-linear partial differential equation}

\author{M. Thamban Nair and  Samprita Das Roy}
\address{Department of Mathemtics, Indian Institute of Technology Madras, Chennai 600036, India.}
\email{mtnair@iitm.ac.in; samprita.dasroy@gmail.com}
\maketitle

\begin{abstract}
We consider a parameter identification problem associated with  a quasi-linear elliptic Neumann boundary value problem involving a parameter function $a(\cdot)$ and the solution $u(\cdot)$, where the problem is to  identify $a(\cdot)$ on an interval $I:= g(\G)$  from the  knowledge of the solution $u(\cdot)$ as $g$ on $\G$, where $\G$ is a given curve on the boundary of the domain $\O\subseteq \R^3$ of the problem and $g$ is a continuous function. The inverse problem is formulated as a problem of solving an operator equation involving a compact operator depending on the data,  and for obtaining stable approximate solutions under noisy data, a new regularization method is considered. The derived error estimates are similar to, and in certain cases better than, the classical Tikhonov regularization considered in the literature in recent past.
\end{abstract}



\tableofcontents

\section{Introduction}\label{sec-1}

Let $\O$ be a bounded domain in $\R^3$ 
with $C^{1,1}$ boundary. Consider the problem of finding a weak solution $u\in H^1(\O)$ of the partial differential equation
\beq\label{eq-1.1}- \nabla . (a(u) \nabla  u) = 0 \q\h{in }\, \O\eeq
with boundary condition
\beq\label{eq-1.2}a(u) \frac{\partial u}{\partial \nu}  = j\q\h{on }\, \partial \O,\eeq
where $a\in H^1(\R)$ and $j\in L^2(\partial \O)$. It is known that such a solution $u$ exists if $a\geq \kappa_0>0$  a.e. for some constant $\kappa_0$  and $\int_{\partial \O} j=0$ (see \cite{Showalter}, \cite{Kugler}). It is also known that, under an additional assumption that $j\in W^{(1-1/p),p}(\partial \O)$ with $p>3$, $u \in C^1(\overline{\O})$ (cf. \cite{Egger}).
One can come across this type of problems in the steady state heat transfer problem with $u$ being the temperature,  $a$ the thermal conductivity which is a function of the temperature, and $j$  the heat flux applied to the surface.

In this paper we consider one of the inverse problems associated with  the above direct problem, namely the following:
\vsq
\noi
{\it {\bf Problem (P):}\, Let $\g  :[0,1] \to \partial \O$ be a $C^1$- curve on $\partial \O$ and $\Gamma$ be its range, that is, $\Gamma:= \gamma([0,1])$. Given $g:\Gamma \to \R$  such that $g\circ \gamma \in C^1([0,1])$ and $j\in W^{(1-1/p),p}(\partial \O)$ with $p>3$ and $\int_{\partial \O} j=0$, the problem is to identify  an $a\in H^1(\R)$ on $I:=g(\G)$ such that the corresponding  $u$ satisfies (\ref{eq-1.1})-(\ref{eq-1.2}) along with the  requirement
\beq\label{eq-1.3} u=g \q \h{on } \Gamma.\eeq }

It is known that, with only the knowledge of $u=g$ on $\Gamma$, the parameter $a$ can be identified uniquely only on $I$\,(cf.\cite{Can}). In the following we shall use the same notation for $a$ for $a\in H^1(\R)$ and its restriction $a_{_I}\in H^1(I)$.

We shall see  that Problem (P) is ill-posed, in the sense that the solution $a_{|_{I}}$ does not depend continuously on the data $g$ and $j$ (see Sections \ref{sec-2}). To obtain a stable approximate solution for  Problem (P),  we use a new regularization method which is different from some of the standard ones in the literature. We discuss this method in Section 3.

The existence and uniqueness of the solution for  Problem (P) is known  under some additional conditions on $\gamma $ and $g$, as specified in Section \ref{sec-2} (c.f. \cite{Can, Egger}). In \cite{Kugler} and \cite{Egger} the problem of finding a stable approximate solution  of the problem is studied by employing Tikhonov regularization with noisy data. In \cite{Kugler}, with the noisy data
$ g^\d$, in place of $g$, satisfying $ \|g-g^\d\|_{L^2(\G)}\leq \d,$  convergence rate $\|a-a^\d\|_{H^1(I)} = O(\sqrt{\delta})$ is obtained whenever  $a\in H^4(I)$ and its trace is Lipschitz on $\partial \O$, where $a^\d$ is the approximate solution obtained via Tikhonov regularization. In \cite{Egger}, the rate $\|a-a^\d\|_{L^2(I)}=O(\sqrt\d)$  is obtained without the additional assumption on $a$. Moreover, here the noisy data $j^{\d}$ belonging to $W^{1-1/p,p}(\partial \O)$ with $p>3$, and satisfying $\|j-j^{\d} \|_{L^2(\partial \O)} \leq \d$, is also considered along with the noisy data $g^{\d}$ as considered in \cite{Kugler}.
It is stated in \cite{Egger} that ``{the rate $O(\sqrt\d)$ is possible with respect to $H^1$-norm, provided some additional smoothness conditions are satisfied}"; however, the details of the analysis is missing.

Under our newly introduced method, we obtain the above type of error estimates using appropriate smoothness assumptions.
In particular we prove that, if $g_1\in \R$ is such that $I=[g_0, g_1]$ and if $a(g_1)$ is known or is approximately known, and the perturbed data
$j^{\d}$ and $g^{\d}$ belong to $W^{1-1/p,p}(\partial \O)$ for $p>3$ and $C^1(\Gamma)$, respectively, satisfying $$\|j-j^{\d} \|_{L^2(\partial \O)} \leq \d,\q \|g-g^{\d}\|_{W^{1,\infty}(\Gamma)} \leq \d,$$
then the convergence rate is  $O(\sqrt{\d})$ with respect to 
$L^2$-norm. With additional assumption that the exact solution is in $H^3(I)$ we obtain a convergence rate  $O({\d}^{2/3})$ with respect to $L^2$-norm.  Again, in particular, if $g\circ\gamma$ is in $H^4([0,1])$, the rate $O({\d}^{2/3})$ \w.r.t. $L^2$-norm  is obtained under  a weaker condition on perturbed data $g^{\d}$, namely,
$g^{\d}\in L^2(\Gamma)$ with $\|g-g^{\d}\|_{L^2(\Gamma)} \leq \d$. Also, in the new method we do not need the assumption on $g^{\d}$ made in \cite{Egger} which is $g^{\d}(\Gamma) \subset g(\Gamma)$.
Thus some of the estimates obtained in this paper are improvements over the known estimates, and are also better than the expected best possible estimate, namely $O(\d^{3/5})$, in the context of Tikhonov regularization, as mentioned in \cite{Egger}.

The paper is organized as follows: In Section \ref{sec-2} we present a theorem which characterize the solution of the inverse Problem (P) in terms of the solution of the Laplace equation with an appropriate Neumann condition. Also,  the inverse problem is represented as the problem of solving a linear operator equation, where the operator is written as a composition of three injective bounded operators one of which is a compact operator, and prove some properties of these operators. The new regularization method is defined in Section \ref{sec-3}, and error estimates with noisy as well as exact data are derived. In Section \ref{sec-4} we present error analysis with some relaxed conditions on the perturbed data. In Section \ref{sec-5} a procedure is described to relax a condition on the exact data and corresponding error estimate is derived. In Section \ref{sec-6} we illustrate the procedure of obtaining a stable approximate solution to Problem (P).

\section{Operator Theoretic Formulation}\label{sec-2}

Throughout the paper we denote by $I$  the range of the function $g:\Gamma \to \R$, and write it as $I=[g_0, g_1]$, that is  $g_0$ and $g_1$ are the left and right end-points of the closed interval $g(\gamma([0,1]))$.

The  following theorem, proved in \cite{Egger}, helps us to identify the solution of Problem (P).

\bt\label{th-existence} Let $j$, $g$ and $\gamma$ be as specified in Problem(P). Then, Problem(P) has a unique solution $a\in H^1(I)$, and it is the  unique $a\in H^1(I)$  such that
\beq \label{eq-2.4} v(\gamma(s))= \int _{g_0} ^{g(\gamma (s))} a(t) dt \q \forall s \in [0,1], \eeq
where $v\in C^1(\overline{\O})$ is the unique function which satisfies
\beq\label{eq-2.1} -\bigtriangleup v=0 \q \h{in } \O, \eeq
 \beq\label{eq-2.2} \frac{\partial v}{\partial \nu}=j \q \h{on } \partial\O \eeq
and
\beq\label{eq-2.3}\int_{\O} v =0.\eeq
\et

It is known that  if  $j\in W^{1-1/p, p}(\partial \O)$ for $p>3$, then $v$ satisfying (\ref{eq-2.1})-(\ref{eq-2.2}) belongs to $W^{2,p}(\O)$, and
\beq \label{ineq-trace-6} \|v^{j}\|_{W^{2,p}(\O)} \leq C \| j\|_{L^2(\partial \O)}\eeq
for some constant $C>0$ (see Theorem 2.4.2.7 and 2.3.3.2 in \cite{Gris}).

In view of Theorem \ref{th-existence},  the inverse Problem (P) can be restated as follows:
Given  $j$ and $g$ as in Problem (P),  let $v \in C^1(\overline{\O})$  be  the function  satisfying (\ref{eq-2.1}), (\ref{eq-2.2}) and (\ref{eq-2.3}). Then, $a\in H^1(I)$ is the solution of Problem(P) if and only if
 $$\int_{g_0}^{g(\g(s))} a(t) dt = v(\g(s)),\q s\in [0, 1].$$
The above equation can be represented as an operator equation
\beq\label{ill-1} Ta = v^j\circ \g,\eeq
where $v^j$ is the solution of (\ref{eq-2.1})-(\ref{eq-2.3}) and the operator  $T: L^2(I)\to L^2[0, 1]$ is defined by
\beq\label{ill-2} (Tw)(s) = \int_{g_0}^{g(\g(s))} w(t) dt, \q w\in L^2(I),\, s\in [0, 1].\eeq

\bt\label{Th-cpt}
The operator $T$ defined in (\ref{ill-2}) is an injective  compact operator of infinite rank.
In particular, $T: H^1(I)\to L^2[0, 1]$ is a compact operator of infinite rank.\et

\bpf

Note that for every $w\in L^2(I)$ and for every $s, \t\in [0, 1]$, we have
$$
|(Tw)(s)- (Tw)(\t)|  \leq   \|w\|_{L^2(I)} |(g\circ\g)(s) - (g\circ\g)(\t)|^{1/2}.$$
Since $g\circ \g$ is continuous,  the set $\{Tw: \|w\|_{L^2(I)}\leq 1\}$ is equicontinuous and uniformly bounded in $C[0,1]$. Hence,  $T$ is a compact operator from $L^2(I)$ to $C[0, 1]$. Since, the inclusion $C[0, 1]\subseteq L^2[0, 1]$ is continuous, it follows that $T: L^2(I)\to L^2[0, 1]$ is also a compact operator. We note that $T$ is injective.
 Hence, $T$  is  of infinite rank.
\epf

It is to be observed that the compact operator $T$ defined in (\ref{ill-2}) depends on $g$. Thus, the problem of solving operator equation (\ref{ill-1}) based on the data $(g, j)$ is  non-linear as well as ill-posed.

In order to consider our  new regularization method for obtaining stable approximate solutions, we  represent the operator $T$ in (\ref{ill-1}) as a composition of three operators, that is,
$$T=T_3T_2T_1,$$
where, for $r\in \{0, 1\}$,
$$T_1:H^r(I) \to H^{r+1}(I),\q  T_2:H^{r+1}(I) \to L^2(I),\q  T_3:L^2(I) \to L^2([0,1])$$
are defined as follows:

\begin{eqnarray}\label{op-T1} T_1(w)(\t) &:=& \int_{g_0} ^\t w(t) dt, \q   w\in H^r(I), \:\: \t \in I,\\
\label{op-T2} T_2(w) & :=& w, \q  w \in H^{r+1}(I),\\
\label{op-T3} T_3(w) &:=& w\circ g \circ \gamma, \q w\in L^2(I) .
\end{eqnarray}
Clearly, $T_1, T_2, T_3$ are linear operators  and
$$(T_3T_2T_1w)(s) = \int_{g_0}^{g(\g(s))} w(t) dt = (Tw)(s),\q s\in [0, 1].$$
Here,  we used the convention that $H^0(I):=L^2(I)$.

By the above representation of $T$, the operator equation (\ref{ill-1}) can be split into three equations:
\begin{eqnarray}
 T_3(\z) &=& v^j\circ \g,\label{eq-split-1}\\
 T_2(b) &=& \z,\label{eq-split-2}\\
 T_1(a) &=& b.\label{eq-split-3}
 \end{eqnarray}
To prove some properties of the operators  $T_1, T_2, T_3$, we specify the requirements  on $j, g$ and $\g$, namely the following.

\begin{Assumption}\label{assum-1}\rm
Let $j\in W^{(1-1/p),p}(\partial \O)$ with $p>3$ and $\int_{\partial \O} j=0$. Let $\g  :[0,1] \to \partial \O$ be a $C^1$- curve on $\partial \O$  and $g:\Gamma \to \R$ be such that $g \in C^1(\Gamma)$,
\beq\label{A-1} C_{\gamma} \leq |\gamma'(s)|\leq C'_{\gamma} \q \forall s\in[0,1],\eeq
\beq\label{A-2} C_g \leq |g'(\gamma(s))|\leq C'_g \q \forall s\in[0,1], \eeq
for some positive constants $C_{\gamma}$, $C'_{\gamma}$, $C_g$ and $C'_g$.
\end{Assumption}

Next we state a result from analysis  which will be used in the next result and also in many other results that follow.

\bl\label{lem-meas}
Let $h_1$ and $h_2$ be two continuous functions on intervals $J_1$ and $J_2$ respectively, such that $h_2(J_2)=J_1$. Also, let $h'_2$ be continuous with $h'_2\neq 0$. Then, we have the following.
$$\int_{J_2} h_1(h_2(x)) dx = \int_{J_1} \frac{h_1(y)}{|h'_2(h_2^{-1}(y))|} dy .$$
\el

We shall also make use of the following proposition.

\bp
Let $C_g,\,  C_\gamma, \, C'_g\,  C'_\gamma$ be as in Assumption \ref{assum-1}. Then for any $w\in L^2(I)$,
\beq \label{T_3-chain-rule} {C_g C_\gamma} \int_0 ^1 |w(g(\gamma(s)))|^2 ds  \leq  \int_I |w(y)|^2 dy \leq
C'_g C'_\gamma  \int_0 ^1 |w(g(\gamma(s)))|^2 ds .\eeq
\ep

\bpf
By Lemma \ref{lem-meas} and the inequalities (\ref{A-1}) and (\ref{A-2}) in Assumption \ref{assum-1}, we have 
$$\int_0 ^1 |w(g(\gamma(s)))|^2 ds
= \int_{g_0}^{g_1} \frac{|w(y)|^2}{|g'(g^{-1}(y)) {\gamma}'({\gamma}^{-1}(g^{-1}(y)))|}dy
                       \leq  \frac{1}{C_g C_\gamma} \int_I |w(y)|^2 dy,$$
$$\int_{g_0}^{g_1} |w(y)|^2 dy
=  \int_0 ^1 |w(g(\gamma(s)))|^2 |g'(\gamma (s)) {\gamma}'(s)| ds
 \leq  C'_g C'_\gamma  \int_0 ^1 |w(g(\gamma(s)))|^2 ds.$$
From the above, we obtain the required inequalities in (\ref{T_3-chain-rule}).
\epf

\bt\label{Th-op-composition}
{Let $r\in \{0, 1\}$, and let
$$T_1:H^r(I) \to H^{r+1}(I),\q T_2:H^{r+1}(I) \to L^2(I),\q  T_3:L^2(I) \to L^2([0,1])$$
be defined as in
(\ref{op-T1}), (\ref{op-T2}) and (\ref{op-T3}), respectively. Then, $T_2$ is a compact operator, and for every $w\in L^2(I)$,
\beq \label{T_1-bdd-below-1}\|w\|_{H^r(I)} \leq \|T_1(w)\|_{H^{r+1}(I)} \leq  (1+ \sqrt{g_1-g_0}) \|w\|_{H^r(I)},\eeq
\beq \label{T_3-bdd-below-2} {C_g C_\gamma} \|T_3(w)\|_{L^2(I)} \leq \|w\|_{L^2(I)} \leq C'_g C'_\gamma \|T_3(w)\|_{L^2([0,1])}, \eeq
In particular, $T_1$ and $T_3$ are bounded operators with bounded inverse from their ranges.}
\et

\bpf Since $H^1(I)$ and $H^2(I)$ are compactly embedded in $L^2(I)$ (cf. \cite{Kes}), $T_2$ is a compact operator of infinite rank. Now, let $w\in H^1(I)$ and $\t \in I$. Then
$$|T_1(w)(\t)|\leq \int_{g_0} ^\t |w(t)| dt \leq \|w\|_{L^2(I)}\sqrt{g_1-g_0},$$
so that
$$\|T_1(w)\|_{L^2(I)}  \leq \|w\|_{L^2(I)} {\sqrt{g_1-g_0}}.$$
Hence, using the fact that $(T_1(w))'=w$ and $(T_1(w))''=w'$,  we have
$$\|w\|_{L^2(I)} \leq   \|T_1(w)\|_{L^2(I)}+\|w\|_{L^2(I)}  \leq (1+ \sqrt{g_1-g_0}) \|w\|_{L^2(I)}$$
so that
$$\|w\|_{L^2(I)} \leq   \|T_1(w)\|_{H^1(I)} \leq (1+ \sqrt{g_1-g_0}) \|w\|_{L^2(I)},$$
$$\|w\|_{H^1(I)} \leq   \|T_1(w)\|_{H^2(I)} \leq (1+ \sqrt{g_1-g_0}) \|w\|_{H^1(I)},$$
Thus, (\ref{T_1-bdd-below-1}) is proved.

By the inequalities in (\ref{T_3-chain-rule}) we obtain
\beq \label{T_3-bdd-below-2} {C_g C_\gamma} \|T_3(w)\|_{L^2([0,1])} \leq \|w\|_{L^2(I)} \leq C'_g C'_\gamma \|T_3(w)\|_{L^2([0,1])} \eeq
for every $w\in L^2(I)$. The inequalities in (\ref{T_1-bdd-below-1}) and (\ref{T_3-bdd-below-2}) also show that $T_1$ and $T_3$ are bounded operator with bounded inverse from their ranges.
\epf

\section{The New Regularization}\label{sec-3}

We know that Problem (P) is ill-posed. We may also recall that the operator equation (\ref{ill-1}) is equivalent to the system of of operator equations (\ref{eq-split-1})-(\ref{eq-split-3}), wherein equation (\ref{eq-split-2}) is ill-posed, since $T_2$ is a compact operator of infinite rank. Thus, in order to regularize (\ref{ill-1}), we shall replace the equation (\ref{eq-split-2}) by a regularized form of it using a family of bounded operators $T_2^\a,\, \a>0$,  which approximates the compact operator $T_2$ in norm.

Note that $T_2:H^2(I) \to L^2(I)$ is defined by
$$T_2(w) = w,\q w\in H^2(I).$$
We consider $T_2^\a$ as a perturbed form of $T_2$, namely,
$T^\a _2 : H^2(I) \to L^2(I)$, defined by
\beq \label{def-reg-T2} T^\a _2 (w)=w-\a w'',\q w\in H^2(I)\eeq
for each $\a>0$.

\bt\label{Th-T2}{For $\a>0$,  let  $T^\a _2 : H^2(I) \to L^2(I)$ be defined as in (\ref{def-reg-T2})
Then,
$$\|T^\a _2 (w)\|_{L^2(I)}\leq  \max\{1, \a\} \|w\|_{H^2(I)},\q  w\in H^2(I).$$
In particular,
$T_2^\a$ is a bounded operator with $\|T_2^\a\|\leq \max\{1, \a\}$.}
Further,
$$\|T_2^\a-T_2\| \to 0\q\h{as}\q \a\to 0.$$
\et

\bpf
We  observe that, for any $w\in H^2(I)$,
$$\|T^\a _2 (w)\|_{L^2(I)}= \|w-\a w''\|_{L^2(I)} \leq \|w\|_{L^2(I)}+ \a \|w''\|_{L^2(I)} \leq \max\{1, \a\} \|w\|_{H^2(I)}.$$
Thus,  $T_2^\a$ is a bounded operator with $\|T_2^\a\|\leq \max\{1, \a\}$ for all $\a>0$.
Further,
$$\|(T^\a _2 -T_2)(w)\|_{L^2(I)} = \|\a w''\|_{L^2(I)} \leq \a \|w\|_{H^2(I)}.$$
Hence, we also have  $\|T_2^\a-T_2\| \to 0$ {as}  $\a\to 0.$
\epf

%


In order to define a regularization family for $T_2$, we introduce the space
\beq
\label{new-space} {\mathcal W}:=\{w\in H^2(I): w(g_0)=0,\,  w'(g_1) = 0\}.\eeq
%
%
Note that, for $w\in H^2(I)$, $w\in {\mathcal W}$ if and only if
$$w(t) = \int_{g_0}^t \xi(s) ds$$
for some $\xi\in H^1(I)$ satisfying $\xi(g_1)=0.$
We prove that $\mathcal{W}$  is a closed subspace of $H^2(I)$ and  $T_2^\a$ as an operator from $\mathcal{W}$ to $L^2(I)$ is bounded below with respect to $H^2(I)$ norm.
%


\bp
The space $\mathcal{W}$ defined in (\ref{new-space}) is a closed subspace of $H^2(I)$ and
$$ (T^\a _2|_{\mathcal W})^* =  Q(T^\a _2)^*,$$
where  $Q: H^2(I)\to H^2(I)$ is the orthogonal projection onto $\mathcal W$.
\ep

\bpf
Let $(w_n)$ in $\mathcal{W}$ be such that $w_n\to w_0$ in $H^2(I)$ for some $w_0\in H^2(I)$.
By a Sobolev imbedding Theorem (cf. \cite{Kes}), $H^2(I)$ is continuously imbedded in the space $C^1(I)$ with $C^1$-norm. Therefore,
$w_0 \in C^1(I)$, and
$$\sup_{t\in I} \{|w_n(t)-w_0(t)| + |w_n'(t)-w'_0(t)|\} \to 0 \q \h{as } n \to \infty.  $$
Also,
$$|w_n(g_0)-w_0(g_0)| \leq \sup_{t\in I} \{|w_n(t)-w_0(t)| + |w'_n(t)-w'_0(t)|\} \q \forall n\in \N $$
and
$$|w'_n(g_1)-w'_0(g_1)| \leq \sup_{t\in I} \{|w_n(t)-w_0(t)| + |w'_n(t)-w'_0(t)|\} \q \forall n\in \N. $$
Thus, since $w_n \in \mathcal{W}$, in particular
$$|w_0(g_0)|= \lim_{n\to \infty} w_n(g_0)=0\q\h{and}\q  |w'_0(g_1)|= \lim_{n\to \infty} w'_n(g_1)=0.$$
Hence $w_0 \in \mathcal{W}$. Thus $\mathcal{W}$ is closed.
Now, let $Q:H^2(I) \to H^2(I)$ be the orthogonal projection onto $\mathcal{W}$. Then,
 for $y\in L^2(I)$ and $w\in \mathcal{W}$ we have,
%
$$ \<Q{(T^{\a}_2)}^*(y),w\>_{H^2(I)}
= \<y, (T^{\a}_2)Qw\>_{L^2(I)}=
                               \<y,(T^{\a}_2|_{\mathcal{W}})w\>_{L^2(I)}
                                = \<{(T^{\a}_2|_{\mathcal{W}})}^*y,w\>_{H^2(I)}
$$

Hence we have $(T^\a _2|_{\mathcal W})^* =  Q(T^\a _2)^*$.
\epf
Let us see some other properties of the space $\mathcal{W}$ which shall be used in order to construct the regularization method.

\bp \label{prop-projection}
 Let $\a >0$.
 Let $L: H^2(I) \to H^2(I)$ be defined by 
$$Lx(t)= x'(g_1)\sqrt{\a}\left[ \frac{e^{(\frac{t-g_0}{\sqrt{\a}})}-e^{-(\frac{t-g_0}{\sqrt{\a}}) }}{e^{(\frac{g_1-g_0}{\sqrt{\a}})} + e^{-(\frac{g_1-g_0}{\sqrt{\a}})}}\right] + x(g_0)\left[ \frac{e^{(\frac{t-g_1}{\sqrt{\a}})}+e^{-(\frac{t-g_1}{\sqrt{\a}}) }}{e^{(\frac{g_0-g_1}{\sqrt{\a}})} + e^{-(\frac{g_0-g_1}{\sqrt{\a}})}}\right]$$
for every $x \in H^2(I)$, $t \in I$. Then we have the following.
\ben
 \i[\rm(i)] For any $x\in H^2(I)$, $Lx \in C^{\infty}(I) \subset H^2(I)$, $\a(Lx)''=Lx$ and  $Lx \in N(T^{\a}_2)$. 
\i[\rm(ii)] $L$ is a bounded linear operator.
 \i[\rm(iii)] The map $id-L$ is a projection onto $\mathcal{W}$, where $id$ is the identity map on $H^2(I)$. 
\een
\ep

\bpf
Clearly, $L$ is a linear operator, and for any $x\in H^2(I)$, we have $Lx \in C^{\infty}(I) \subset H^2(I)$ and $\a (Lx)''=Lx$.
To show that  $L$ is continuous, let  $(x_n)$ be a sequence in $H^2(I)$ such that $\|x_n -x\|_{H^2(I)} \to 0$ for some $x\in H^2(I)$. By a Sobolev imbedding Theorem (cf. \cite{Kes}), $H^2(I)$ is continuously imbedded in the space $C^1(I)$ with $C^1$-norm, and so we have $|x_n(g_0)-x(g_0)| \to 0$ and $|x'_n(g_1)-x'(g_1)| \to 0$ as $n \to \infty$. Using this, it can be shown that $L$ is continuous. 
Now again by definition of $L$, for any $x\in H^2(I)$ we have 
\beqarray
(x-Lx)(g_0) &=& x(g_0)- Lx(g_0)=x(g_0)-x(g_0)=0,\\
(x-Lx)'(g_1)&=& x'(g_1)- (Lx)'(g_1)= x'(g_1)-x'(g_1)=0,
\eeqarray
so that 
$(id-L)(x-Lx)=x-Lx- L(x-Lx)=x-Lx.$
Hence, using the definition of the space $\mathcal{W}$, we have $id-L$ is a projection onto $\mathcal{W}$. 
\epf

We shall use the notation  
\beq\label{C_L} C_L:=\|id-L\|,\eeq
where $L$ is the bounded operator as in Proposition \ref{prop-projection}.

\bt\label{reg-bdd-below-1}
Let $0<\a<1$. Then, for every $w\in {\mathcal W}$,
\beq\label{bdd-below-1}\|T^\a _2(w)\|_{L^2(I)}\geq \a \|w\|_{H^2(I)},\eeq
\beq\label{bdd-below-2}\|T^\a _2(w)\|_{L^2(I)}\geq \sqrt\a \|w\|_{H^1(I)}.\eeq
%
\et

\bpf
First we observe, by integration by parts, that for $w_1,w_2 \in {\mathcal{W}}$,
$\int_I w_{1}w''_{2}= -\int_I w'_{1}w'_{2}.$
Hence, for every $w\in {\mathcal{W}}$,
\beqarray \|T^\a _2(w)\|^2_{L^2(I)} &=& \int_{g_0} ^{g_1} |w-{\a} w''|^2 \\
                                  &=& \int_{g_0} ^{g_1} |w|^2
  + {\a}^2 \int_{g_0} ^{g_1} |w''|^2 - 2{\a} \int_{g_0} ^{g_1} ww''                                \\
  &=& \int_{g_0} ^{g_1} |w|^2   + {\a}^2 \int_{g_0} ^{g_1} |w''|^2 + 2{\a} \int_{g_0} ^{g_1} |w'|^2.
\eeqarray
Since   $0<\a<1$, for every $w\in \mathcal{W}$,
$$\int_{g_0} ^{g_1} |w|^2+ {\a}^2 \int_{g_0} ^{g_1} |w''|^2 + 2{\a} \int_{g_0} ^{g_1} |w'|^2\geq
 \a^2 \|w\|_{H^2(I)}^2,$$
$$\int_{g_0} ^{g_1} |w|^2+ {\a}^2 \int_{g_0} ^{g_1} |w''|^2 + 2{\a} \int_{g_0} ^{g_1} |w'|^2\geq
 \a \|w\|_{H^1(I)}^2.$$
 This completes the proof.
\epf

At this point let us note that,  by (\ref{bdd-below-1}),  $T_2^\a$ is  is bounded below on $\mathcal{W}$.
Henceforth, we shall use the same notation for  $T_2^\a$  and its restriction to $\mathcal{W}$, that is,
\beq \label{T2-rest} T^\a _2 (w)=w-\a w'',\q w\in \mathcal{W}\eeq
and the adjoint of this operator will be denoted $(T^\a _2)^*$.

In the following, we  use the notations  $R(S)$ and $N(S)$ for the range and null space, respectively,  of the operator $S$.
\bl \label{lem-bdd-bel}
Let $H_1$ and $H_2$ be Hilbert spaces  and let  $S: H_1\to H_2$ be a bounded linear operator  with  closed range. Then,
\beq\label{prop-closedrange} R(S^*S) = R(S^*)\eeq
Suppose,  in addition, that there exist $c>0$  such that
$\|Sx\|\geq c\|x\|$  for all\, $x\in H_1.$
Then
\beq\label{prop-closedrange-1}
\|S^*Sx\|\geq c^2 \|x\|\q\forall\, x\in H_1,\eeq
Further, if $\|\cdot\|_0$ is any  norm on $H_1$ and if  $c_0>0$ is such that
$\|Sx\| \geq c_0 \|x\|_0$  for all\, $x\in H_1$,
then
\beq\label{prop-closedrange-3}\|S^\dagger\,y\|_0\leq  \frac{1}{c_0} \|y\|\q\forall\, y \in H_2,\eeq
where $S^\dagger:= (S^*S)^{-1}S^*$, the generalized inverse of $S$.
\el

\bpf
Clearly,  $R(S^*S) \subseteq  R(S^*)$. Now, let $x\in R(S^*)$, and let $y\in H_2$ be such that $x=S^*y$. Let $y_1\in N(S^*)$ and $y_2\in N(S^*)^\perp$ be such that $y= y_1+y_2$. Hence, $x=S^*y_2$. Since $R(S)$ is closed, $N(S^*)^\perp= R(S)$. Hence, there exists $x_2\in H_1$ such that $y_2=Sx_2$. So, $x=S^*Sx_2\in R(S^*S)$. Thus, $R(S^*) \subseteq  R(S^*S)$. Thus, we have proved (\ref{prop-closedrange}).

Next, suppose that there exist $c>0$  such that
$\|Sx\|\geq c\|x\|$  for all\, $x\in H_1.$
Then for every $x\in H_1$,
$$\|S^*Sx\|\, \|x\| \geq \<S^*Sx, x\>_{H_1} = \|Sx\|^2\geq c^2\|x\|^2.$$
Thus, we obtain (\ref{prop-closedrange-1}).

By (\ref{prop-closedrange-1}), $R(S^*S)$ is closed and $S^*S$ has a bounded inverse from its range and hence, by (\ref{prop-closedrange}),  $(S^*S)^{-1}S^*$ is well defined as a bounded operator from $H_2$ to $H_1$. Since $R(S)$ is closed, it is known that for every $y\in H_2$, there exists $x\in H_1$ such that
\beq\label{prop-closedrange-4}(S^*S)x=S^*y\q\h{and}\q Sx=Py,\eeq
where $P:H_2\to H_2$ is the orthogonal projection onto $\overline{R(S)} = R(S)$, and this $x$ is unique since $S$ and $S^*S$ are bounded below  (see, e.g. \cite{nair-fa}).
Now,  assume that  $\|\cdot\|_0$ is any norm on $H_1$ such that
$\|Sx\| \geq c_0 \|x\|_0$  for all\, $x\in H_1$ for some $c_0>0$.
For   $y\in H_2$, if $x$ is as in (\ref{prop-closedrange-4}), then
$$\|(S^*S)^{-1}S^*y\|_0 = \|x\|_0 \leq \frac{1}{c_0}\|Sx\| = \frac{1}{c_0}\|Py\|\leq \frac{1}{c_0} \|y\|.$$
Thus, we obtain (\ref{prop-closedrange-3}).
\epf

\bcor\label{cor-bdd-bel}
Let  $0<\a<1$ and  $T_2^\a$ be as in (\ref{T2-rest}). Then for every $y\in L^2(I)$,
\begin{eqnarray}
\|{({(T_2^\a)}^*T_2^\a)}^{-1} (T_2^\a)^*y\|_{H^2(I)}
&\leq & \frac{1}{\a} \|y\|_{L^2(I)},\label{ineq-reg-1}\\
\|{({(T_2^\a)}^*T_2^\a)}^{-1} (T_2^\a)^*y\|_{H^1(I)} & \leq &\frac{1}{ \sqrt{\a}} \|y\|_{L^2(I)},\label{ineq-reg-2}
\end{eqnarray}
\ecor

%


\bpf
Taking $H_1={\mathcal W}$ and $H_2=L^2(I)$ in Lemma \ref{lem-bdd-bel}, the   inequalities  in  (\ref{ineq-reg-1}) and (\ref{ineq-reg-2})  follow from (\ref{prop-closedrange-3}) by taking the norm $\|\cdot\|_0$ as  $\|\cdot\|_{H^2(I)}$ and $\|\cdot\|_{H^1(I)}$ respectively, on ${\mathcal W}$ and by  using (\ref{bdd-below-1}) and  (\ref{bdd-below-2}), respectively.
\epf

Let $R_\a: L^2(I)\to {\mathcal W}$ for $\a>0$ be defined by
\beq\label{reg-fly}
R_\a :={({(T^\a _2)}^* (T^\a _2))}^{-1} {(T^\a _2)}^*,\q\a>0.\eeq
We note that, by Corollary \ref{cor-bdd-bel}, $R_\a$ is a bounded operator  from $L^2(I)$ to ${\mathcal W}$ (with respect to the norm $\|\cdot\|_{H^2(I)}$), for each $\a>0$. Since
$(T_2-  T^\a _2)(w) = \a w''$, we have
\beq\label{reg-fly-1}R_\a T_2w -w
= \a R_\a (w'').\eeq

Next, we prove that $\{R_\a\}_{\a>0}$, defined as in (\ref{reg-fly}),   is a regularization family for $T_2: {\mathcal W}\to L^2(I)$. Towards this aim, we first prove the following theorem.

\bt\label{th-reg}
For $\a>0$, let $R_\a$ be as in (\ref{reg-fly}), and let $C_L$ be as in (\ref{C_L}).  Then the following results hold.
\begin{itemize}
\item[(i)] $\|R_\a T_2w \|_{H^2(I)} \leq 2 \|w\|_{H^2(I)}$ for all $w\in \mathcal{W}$.
\item[(ii)] $\|R_\a T_2 w -w \|_{H^2(I)}\leq (1+C_{L})\a\|w''\|_{H^2(I)}$ for all $w$ in $\mathcal{W} \cap H^4(I)$.
\item[(iii)] $\|R_\a T_2 w -w \|_{H^1(I)}\leq \sqrt{\a} \|w''\|_{L^2(I)}$ for all $w$ in $\mathcal{W}$.
\end{itemize}
\et

\bpf
(i)\, Let $w \in \mathcal{W}$. By (\ref{reg-fly-1}), we have
$$\|R_\a T_2w\|_{H^2(I)} = \|w - [w- R_\a T_2(w)]\|_{H^2(I)} = \|w + \a R_\a (w'')\|_{H^2(I)}.$$
Hence, using (\ref{ineq-reg-1}),
$$\|R_\a T_2w\|_{H^2(I)} \leq \|w\|_{H^2(I)}+\a \|R_\a (w'')\|_{H^2(I)}  \leq  \|w\|_{H^2(I)}+ \|w''\|_{L^2(I)}.$$
Thus, $\|R_\a T_2w\|_{H^2(I)} \leq 2  \|w\|_{H^2(I)}$ for every $w \in \mathcal{W}$.

\vsq

(ii)\, Let $w \in \mathcal{W} \cap H^4(I)$. Let us note that $w''$ is in the domain of $T_2$ and hence is in $H^2(I)$ (may not be in $\mathcal{W}$). By Proposition \ref{prop-projection}, $w''-Lw'' \in \mathcal{W}$ and $Lw'' \in N(T^{\a}_2)$.
Thus, using the above fact, along with the fact that $w''$ is in the domain of $T_2$, by  (\ref{reg-fly-1}) and  (i) above, we have
\beqarray
\|R_\a T_2w-w\|_{H^2(I)} &= &  \a \|R_\a (w'')\|_{H^2(I)}\\
                         & = & \a \|R_\a T_2 (w'')\|_{H^2(I)} \\
                         &=& \a \|R_\a [T^{\a}_2 (w'') + \a w'''']\|_{H^2(I)}\\
                         & \leq & {\a}^2 \|R_\a (w'''')\|_{H^2(I)} + \a \|R_\a T^{\a}_2 (w'')\|_{H^2(I)}\\
                         &=& {\a}^2 \|R_\a (w'''')\|_{H^2(I)} + \a \|R_\a T^{\a}_2 (L w'') + R_\a T^{\a}_2[(id-L) (w'')]\|_{H^2(I)}\\
                         & = & {\a}^2 \|R_\a (w'''')\|_{H^2(I)} +\a \|R_\a T^{\a}_2 [(id-L) (w'')]\|_{H^2(I)}\\
                         &=& {\a}^2 \|R_\a (w'''')\|_{H^2(I)} +\a \|(id-L) (w'')\|_{H^2(I)}\\
                         & \leq & \a [\|w''''\|_{L^2(I)} + \|(id-L) (w'')\|_{H^2(I)}].
\eeqarray
Now, since $\|w''''\|_{L^2(I)}\leq  \|w''\|_{H^2(I)}$ and 
$\|(id-L) (w'')\|_{H^2(I)}\leq C_L\|w''\|_{H^2(I)}$, we obtain the required inequality. 

(iii)\, For $w\in {\mathcal W}$,  using (\ref{ineq-reg-2}), we have
$\|R_\a T_2w-w\|_{H^1(I)} =  \a \|R_\a (w'')\|_{H^1(I)} \leq  \sqrt{\a} \|w''\|_{L^2(I)}.$
Thus, the proof is complete. 
\end{proof}

%


\bl \label{lem-density}
The space $\mathcal{W}\cap H^4(I)$ is dense in $\mathcal{W}$.
\el
\bpf
Let $w\in {\mathcal{W}}$. Since $H^4(I)$ is dense in $H^2(I)$ as a subspace of $H^2(I)$ (cf. \cite{Kes}), there exists  a sequence $(w_n)$ in $H^4(I)$ such that
\beq \label{eq-density-1} \|w_n-w\|_{H^2(I)} \to 0 \q \h{as } n\to 0. \eeq
Now, define $P:H^2(I) \to \mathcal{W}$ by
$$P(w)(t)=w(t)-w(g_0)-w'(g_1)(t-g_0), \q w\in H^2(I) \q \h{and } t\in I. $$
Since  $H^2(I)$ is continuously imbedded in $C^1(I)$ (cf. \cite{Kes}), (\ref{eq-density-1}) implies that
$|w_n(g_0)-w(g_0)| \to 0$  and $ |w'_n(g_1)-w'(g_1)| \to 0$
as $n \to 0.$
Thus, as $I$ is bounded we have
\beq \label{eq-density-2} \|P(w_n)-P(w)\|_{H^2(I)} \to 0 \q \h{as } n\to 0. \eeq
Again by definition of $P$ and $\mathcal{W}$ we have
$Pw_n \in \mathcal{W}\cap H^4(I)$ and $Pw=w$. Hence from (\ref{eq-density-1}) and (\ref{eq-density-2}) we have the proof.
\epf


\bt \label{th-reg-conv}
Let $w\in \mathcal{W}$, and let $\{R_{\a}\}_{\a>0}$ be as in (\ref{reg-fly}). Then
 $$\|R_\a T_2 w - w\|_{H^2(I)} \to 0 \q \h{as } \a \to 0.$$
 In particular, $\{R_{\a}\}_{\a>0}$ is a regularization family for $T_2$.
\et
\bpf
By Theorem \ref{th-reg}, $(R_\a T_2)$ is a uniformly bounded family of operators from ${\mathcal W}$ to ${\mathcal W}$ and
$\|R_\a T_2 w - w\|_{H^2(I)} \to 0$ as $\a\to 0$ for every $x\in \mathcal{W}\cap H^4(I)$.
Since   $\mathcal{W}\cap H^4(I)$ is dense in ${\mathcal W}$  (see  Lemma \ref{lem-density}), by a result in functional analysis (see Theorem 3.11 in \cite{nair-fa}), we obtain 
$\|R_\a T_2 w - w\|_{H^2(I)} \to 0$ as $\a \to 0$
for every $w\in {\mathcal W}$. Thus $\{R_{\a}\}_{\a>0}$ is a regularization family for $T_2$.
\epf

Throughout, we assume that $a_0\in H^1(I)$ is the unique solution of  Problem (P).
Thus, equations (\ref{eq-split-1})-(\ref{eq-split-3}) have solutions namely, $\zeta_0, b_0$ and $a_0$, respectively. That is,
\begin{eqnarray}
 T_3(\z_0) &=& v^j\circ \g,\label{eq-split-1-0}\\
 T_2(b_0) &=& \z_0,\label{eq-split-2-0}\\
 T_1(a_0) &=& b_0.\label{eq-split-3-0}
 \end{eqnarray}
Having obtained the regularization family $\{R_\a\}_{\a>0}$ for $T_2$ as in (\ref{reg-fly}), we may replace the solution $b_0$ of the equation $(\ref{eq-split-2})$  by
$$b_\a:= R_\a \zeta_0.$$
Thus, we may define the regularized solution $a_\a$ for Problem (P) as the solution of (\ref{eq-split-3}) with $b_0$ replaced by $b_\a$. Thus the regularized solution $a_\a$ for Problem (P) is defined along the following lines:
\begin{eqnarray}
 T_3 (\zeta_0) &=& v^j \circ \g,\label{eq-split-4}\\
 (T_2^\a)^*T_2^\a (b_\a ) &=& (T_2^\a)^*\zeta_0,\label{eq-split-5}\\
 T_1(a_\a) &=& b_\a.\label{eq-split-6}
 \end{eqnarray}
Since $b_\a \in \mathcal{W} \subset R(T_1)$, each of the  above equations has unique solution. In fact
$\zeta_0=T_2b_0$ with  $b_0=T_1a_0$, where $a_0$ is the unique solution of (\ref{ill-1}). Note that, the operator equation  (\ref{eq-split-5}) has a unique solution because $T_2^\a$ is bounded below, and  (\ref{eq-split-6}) has a unique solution as  $T_1$ is injective with range ${\mathcal W}$, and $b_\a\in {\mathcal W}$. Hence we have,
$a_{\a}(g_1)=0$. Thus to obtain convergence of $\{a_{\a}\}$ to $a_0$ as $\a \to 0$, it is necessary that $a_0(g_1)=0$. Therefore, in this section, we assume that,
\beq \label{A-4} a_0(g_1)=0.\eeq
 We shall relax this condition in Section \ref{sec-5}, by appropriately redefining regularized solutions.

\subsection{Error estimates under exact data}

For $\a>0$, let $a_\a$ be defined via equations (\ref{eq-split-4})-(\ref{eq-split-6}). Also, let $a_0$ be the unique solution to Problem(P) satisfying (\ref{A-4}).
Then, we look at the estimates for the error term $(a_0-a_\a)$ in both $L^2(I)$ and $H^1(I)$ norms in the following theorem.

\bt 
The following results hold.
\ben
\i $\|a_0-a_\a\|_{H^1(I)}\to 0$ {as} $\a\to 0.$
\i $\|a_0- a_\a\|_{L^2(I)}\leq  \sqrt \a \|a'_0\|_{L^2(I)}$.
\i If $a_0\in H^3(I)$, then with $C_L$ is as in (\ref{C_L}), $\|a_0-a_\a\|_{H^1(I)}\leq (1+ C_L) \a \|a_0'\|_{H^2(I)}$.
\een
\et

\bpf
By our assumption,  $a_0(g_1)=0$. Therefore, by definition of $T_1$ and the space $\mathcal{W}$, we have

$b_0=T_1(a_0)\in \mathcal{W}.$
Now let us first observe that, by the definition of $b_\a$
$$ T_1(a_0)-T_1(a_\a) =  b_0-b_{\a} = b_0- R_\a \z_0  = b_0- R_\a T_2 b_0.$$
Hence, by the inequality (\ref{T_1-bdd-below-1}), for $r\in \{0, 1\}$,
we have,
\beq\label{error-1} \|a_0-a_{\a}\|_{H^r(I)} \leq \|T_1(a_0)-T_1(a_\a)\|_{H^{r+1}(I)} =
\|b_0- R_\a T_2 b_0\|_{H^{r+1}(I)},\eeq
and hence, by Theorem \ref{th-reg-conv}, $\|a_0-a_{\a}\|_{H^1(I)} \to 0$ as $\a\to 0$.  Thus we have proved  (1).

Also, since $b_0\in {\mathcal W}$, from  (\ref{error-1}) and  Theorem \ref{th-reg}(iii), we have
$$\|a_0-a_{\a}\|_{L^2(I)} \leq   \|T_1(a_0)-T_1(a_\a)\|_{H^1(I)}
= \|b_0- R_\a T_2 b_0\|_{H^1(I)}\leq \sqrt{\a}\|b_0''\|_{L^2(I)} = \sqrt{\a}\|a_0'\|_{L^2(I)}.$$
which proves (2). Now, let $a_0\in H^3(I)$. Then $b_0\in H^4(I)$. Since $b_0\in \mathcal{W}$, we have $b_0\in \mathcal{W} \cap H^4(I)$.
 Hence proof of (3) follows from (\ref{error-1}) and Theorem \ref{th-reg} (ii).
\epf

%
%
%
%
%
%
%
%

\subsection{Error estimates under noisy data}
In practical situations the observations of the data $j$ and $g$ may not be known accurately and we may have some noisy data  instead. In this section we assume that the noisy data $g^{\e}$ and $j^{\d}$  are such that
\beq \label{pert-0} g^{\e} \in C^1(\Gamma),\q  \q j^{\d} \in W^{1-1/p,p}(\partial \O),\, p>3 \eeq  satisfying
\beq \label{pert-1} \|g-g^\e \|_{W^{1,\infty}(\Gamma)} \leq \e, \eeq
\beq \label{pert-2} \|j-j^\delta \|_{L^2(\partial \O)} \leq \delta \eeq
for some known noise level $\e$ and $\d$, respectively.
At this point let us note that a weaker condition on perturbed data $j^{\d}$, for example $j^\d \in L^2(\partial \O)$, is not very feasible to work with, in this problem. This is because, in that case the corresponding solution $v^{j^{\d}}$ of (\ref{eq-2.2})-(\ref{eq-2.3}) with $j^{\d}$ in place of $j$, is not continuous and hence its restriction on $\Gamma$ does not make sense. In practical situations if such a perturbed data arise we may work with its appropriate approximation which is in $W^{1-1/p,p}(\partial \O)$ with $p>3$. For the perturbed data $g^\e$, in the next section we consider the case when it is in a more general space which is $L^2(\Gamma)$.

Corresponding to the data $j, j^\d$ as above, we denote
\beq\label{data-f} f^j:=v^j\circ \g,\q f^{j^\d}:= v^{j^\d}\circ \g.\eeq

\bl\label{lem-trace-0}
Let $\gamma_0$ be a $C^1$ curve on $\R^2$ and let $\Gamma_0 = \{(x, \gamma_0(x)) \in \R^2 : d_0 \leq x\leq d_1 \}$
for some $d_0, d_1$ in $\R$ with $d_0<d_1$. Then
\beq \label{ineq-trace-1}\|w\|_{L^2(\Gamma_0)} \leq  \|w\|_{H^1(\R^2)}, \q \forall w \in H^1(\R^2).\eeq
\el
\bpf
Let $w\in C^{\infty}_c(\R^2)$. Then, using H\"{o}lder's inequality we have
\beqarray
\|w\|^2_{L^2(\Gamma_0)} & = & \int_{\Gamma_0} (w(z))^2 dz = \int_{d_0}^{d_1} (w(x,\gamma_0(x)))^2 dx \\
                            &=& \int_{d_0}^{d_1} \left[ -\left( \int_{\gamma_0 (x)}^{\infty } \frac{\partial}{\partial t} (w(x,t))                     ^2 dt \right) \right]  dx  \\
&=&   \int_{d_0}^{d_1} \left( \int_{\gamma_0 (x)}^{\infty} (-2w(x,t) \frac{\partial}{\partial t}w(x,t)) dt \right) dx \\
& \leq &  \int_{d_0}^{d_1}  \left( \int_{\gamma_0 (x)}^{\infty} |w(x,t)|^2 dt  + \int_{\gamma_0 (x)}^{\infty} |\frac{\partial}{\partial t} (w(x,t)|^2 dt \right) dx \\
& \leq & \|w\|^2_{L^2(\R^2)} +  \| \nabla w \|^2_{L^2(\R^2)} \\
                           & \leq &  \|w\|^2_{H^1(\R^2)}.
\eeqarray
Hence, $C^{\infty}_c(\R^2)$ being dense in $H^1(\R^2)$, we have the proof.
\epf

\bl\label{lem-trace-1}
Let $w\in H^{1}(\partial \O)$ and $\g$ be a curve on $\partial \O$ such that $|\g'(t)|$ is bounded away from $0$ as in (\ref{A-1}). Then there exists $C_0>0$ such that
$$\|w\circ \gamma \|_{L^2([0,1])} \leq  C_0  \|w\|_{H^1(\partial \O)}.$$
\el

\bpf
Let $w\in H^{1}(\partial \O)$. Since $\O$ is with $C^1$ boundary,
\beq \label{trace-norm}\|w\|_{H^{1}(\partial \O)} := \sum_{i=1}^m \|\o_i\|_{H^{1}(\R^2)}\eeq for some elements
$\o_1$, $\o_2$,$\cdots$, $\o_m$  $\in H^1(\R^2)$ (cf. \cite{Gris}, \cite{Kes}).
Also, there exists a set $\{\sigma_1,\cdots, \sigma_m\}$ of diffeomorphisms from some neighbourhoods in $\partial \O$ to $\R^2$, which satisfies
\beq \label{trace-norm-1}\|w\circ \gamma \|_{L^2([0,1])}= \sum_{i=1}^m \|\o_i \circ \sigma_i \circ \gamma\|_{L^2([0,1])}.\eeq
For any $i\in \{1,\cdots,m\}$, since $\sigma_i$ is a diffeomorphism  $\sigma_i \circ \gamma $ is a curve in $\R^2$.
Since $|\g'|$ is bounded away from $0$, there exists constant $C_{\gamma} >0$ such that
$|\g'(t)|\geq C_{\g}$ for all $t\in [0,1]$. Also, as $\sigma'_i (\Gamma)$ is compact and $\sigma_i$ is one-one there exists constant $C_{\sigma}>0$ such that $ |\sigma'_i(x)| \geq C_{\sigma}$ for all $x\in \g([0,1])$ and $1\leq i \leq m$. Hence, by Lemma \ref{lem-meas} and (\ref{trace-norm-1}), we obtain
$$\|w\circ \gamma \|_{L^2([0,1])}= \sum_{i=1}^m \|\o_i \circ \sigma_i \circ \gamma\|_{L^2([0,1])} \leq \frac{1}{\sqrt{C_{\gamma}C_{\sigma}}} \sum_{i=1}^m \|\o_i\|_{L^2(\sigma_i(\Gamma))}. $$
Hence, using (\ref{ineq-trace-1}) and (\ref{trace-norm}), we get
$$\|w\circ \gamma \|_{L^2([0,1])} \leq \frac{1}{\sqrt{C_{\sigma} C_{\gamma}}} \sum_{i=1}^m \|\o_i\|_{H^1(\R^2)} = \frac{1}{\sqrt{C_{\sigma} C_{\gamma}}}\|w\|_{H^1(\partial \O)}.$$
This completes the proof.
\epf

\bp \label{p-trace}
Let $\tilde j \in W^{1-1/p,p}(\partial \O)$. Let $v^{\tilde{j}} \in W^{1,p}(\O)$ be the solution of (\ref{eq-2.2})-(\ref{eq-2.3}) with $\tilde{j}$ in place of $j$, such that it satisfies (\ref{eq-2.4}). Then there exists $\tilde C_\g>0$ such that
$$\|v^{\tilde{j}} \circ \gamma \|_{L^2([0,1])} \leq \tilde{C_{\gamma}}\|\tilde{j}\|_{L^2(\partial \O)}.$$
\ep

\bpf
Since $\tilde j$ is in $W^{1-1/p, p}(\partial \O)$, we know that  $v^{\tilde j}\in  W^{2,p}(\O)$ (cf. \cite{Gris}) and
\beq \label{ineq-trace-2} \|v^{\tilde{j}}\|_{W^{2,p}(\O)} \leq C_5 \|\tilde j\|_{L^2(\partial \O)}\eeq
for some constant $C_5>0$ (see inequality \ref{ineq-trace-6}).
By {trace theorem} for Sobolev Spaces (cf. \cite{Gris}), and by continuous imbedding of $W^{(2-1/p),p}(\partial \O)$  into $W^{1,p}(\partial \O)$, we have
$v^{\tilde{j}}|_{\partial \O} \in W^{2-{1}/{p},p}(\partial \O) \subseteq W^{1,p}(\partial \O)$ and
\beq \label{ineq-trace-3}\|v^{\tilde{j}}_{|_{\partial \O}} \|_{W^{1,p}(\partial \O)} \leq C_6\|v^{\tilde{j}}|_{\partial \O} \|_{W^{2-1/p,p}(\partial \O)} \leq C_7\|v^{\tilde{j}}\|_{W^{2,p}(\O)} \eeq
for some constants $C_6, C_7>0$.

Since $p>3$, we have $v^{\tilde{j}}_{|_{\partial\O}} \in H^1(\partial\O)$ and, there exists constant $C_8>0$ such that
$$\|v^{\tilde{j}}_{|_{\partial\O}}\|_{H^1(\partial\O)} \leq C_8\|v^{\tilde{j}}_{|_{\partial\O}}\|_{W^{1,p}(\partial\O)} .$$
Thus, using (\ref{ineq-trace-2}), (\ref{ineq-trace-3}) and with $v^{\tilde{j}}_{|_{\partial \O}}$ in place of $w$ in Lemma \ref{lem-trace-1} we have,
$$\|v^{\tilde{j}}\circ \gamma \|_{L^2([0,1])} \leq \frac{1}{\sqrt{C_{\sigma} C_{\gamma}}}\|v^{\tilde{j}}_{|_{\partial \O}}\|_{H^1(\partial \O)} \leq \frac{C_8}{\sqrt{C_{\sigma} C_{\gamma}}}\|v^{\tilde{j}}_{|_{\partial \O}}\|_{W^{1,p}(\partial \O)} \leq \tilde C_{\gamma}\|\tilde{j}\|_{L^2(\partial \O)},$$
 where $\tilde{C_{\gamma}}= {C_8C_7C_5}/{\sqrt{C_{\sigma} C_{\gamma}}}$.
\epf

\bcor
Let $j$ be as in Assumption \ref{assum-1} and $j^\d$ satisfy (\ref{pert-0}) and (\ref{pert-2}). Let $f$ and $f^{j^\d}$ be as in (\ref{data-f}). Then
\beq \label{pert-3} \|f^j - f^{j^\d} \|_{L^2([0,1])} \leq \tilde{C_{\gamma}} \d,\eeq
where $\tilde{C_{\gamma}}>0$ is as in Proposition \ref{p-trace}.
\ecor

\bpf
By Proposition \ref{p-trace} we have
 $$\|f^j - f^{j^\d} \|_{L^2([0,1])} \leq \tilde{C_{\gamma}} \|j- j^{\d}\|_{L^2(\partial\O)} \leq \tilde{C_{\gamma}}\d.$$
 Hence,
$ \|f^j - f^{j^\d} \|_{L^2([0,1])} \leq \tilde{C_{\gamma}} \d.$
\epf

\bl
For $\e>0$,
\beq \label{ineq-pert-0}C_g-\e\leq |{g^{\e}}'(\gamma(s))| \leq C'_g +\e,\eeq
where $C_g$ and $C'_g$ are as in (\ref{A-2}).
In particular, if \,  $0 < \e \leq {C_g}/{2}$ then
\beq \label{ineq-pert} \frac{C_g}{2} \leq |{g^{\e}}'(\gamma(s))| \leq 2C'_g \, \forall s \in [0,1].\eeq
\el

\bpf
For any $s$ in $[0,1]$, we have
 $$|g'(\gamma(s))|-|g'(\gamma(s))-{g^{\e}}'(\gamma(s))| \leq |{g^{\e}}'(\gamma(s))| \leq  |{g^{\e}}'(\gamma(s)) - g'(\gamma(s))| + |g'(\gamma(s))|.$$
Since
$|g'(\gamma(s))-{g^{\e}}'(\gamma(s))| \leq \|g-g^{\e}\|_{W^{1,\infty}(\Gamma)}<\e,$
by   (\ref{A-2}), we obtain (\ref{ineq-pert-0}). The relations in (\ref{ineq-pert}) are obvious by the assumption on $\e$.
\epf

\brem\label{rem-intervals-1}
Since, $\gamma'$ satisfies (\ref{A-1}), and, $(g^{\e})'$ satisfies (\ref{ineq-pert}) for $\e < {C_g}/{2}$,
$g^{\e}(\Gamma)$ is a non-degenerate closed interval, that is,
$I_{\e}:= g^{\e}(\Gamma) = [g^{\e}_0,\, g^{\e}_1]$  for some $g^{\e}_0, g^{\e}_1$  with $g^{\e}_0 < g^{\e}_1.$
\erem

The following lemma  will help us in showing that $I\cap I_\e$ is a closed and bounded (non-degenerate) interval.

\bl \label{lem-interval-1}
Let $\phi_1$, $\phi_2$ be in $C([\xi_1, \xi_2])$ for some $\xi_1$ and $\xi_2$ in $\R$, and let $\eta>0$ be such that
\beq \label{eq-interval-1} \|\phi_1-\phi_2\|_{L^{\infty}([\xi_1, \xi_2])} \leq \eta. \eeq
Let $I_1 :=\phi_1([\xi_1, \xi_2])=[a_1, b_1]$ and $I_2 := \phi_2([\xi_1, \xi_2])=[a_2, b_2]$ for some $a_1$, $b_1$, $a_2$ and $b_2$ in $\R$.  We assume that $I_1$ and $I_2$ are non-degenerate intervals, that is, $a_1< b_1$ and $a_2< b_2$,
and
\beq \label{A-interval} 2\eta < \min\{(b_1-a_1), (b_2-a_2)\}.\eeq
Then
\beq \label{ineq-interval-1} \max\{|a_1-a_2|,\,  |b_1-b_2|\} \leq \eta \eeq
and $ I_1 \cap I_2=[a, b]$ is a non-degenerate interval, that is, $a <b$.
\el

\bpf
Suppose $a_1< b_1$ and $a_2< b_2$.  Since  $a_1=\phi_1(s_1), a_2=\phi_2(s_2), b_1=\phi_1(s_1'), b_2=\phi_2(s_2')$,    for some $s_1, s_2, s_1', s_2' \in [\xi_1, \xi_2]$, and since
$a_1\leq \phi_1(s_2), a_2\leq \phi_2(s_1), b_1\geq \phi_1(s_2')$ and $b_2\geq \phi_2(s_1')$,
we obtain
\beq \label{A-interval-1}|a_1-a_2| \leq \|\phi_1-\phi_2\|_{L^{\infty}([\xi_1, \xi_2])} \leq \eta, \eeq
\beq \label{A-interval-2}|b_1-b_2| \leq \|\phi_1-\phi_2\|_{L^{\infty}([\xi_1, \xi_2])} \leq \eta. \eeq
Thus, (\ref{ineq-interval-1}) is proved.

To prove the remaining, let us first consider the case $a_1 \leq a_2$.  Then, $I_1\cap I_2=[a_2,\tilde b]$, where $\tilde b:=\min\{b_2, b_1\}$. Note that,  by (\ref{A-interval}) and (\ref{A-interval-1}), we have
$$b_1-a_2 = (b_1-a_1) -(a_2-a_1) \geq 2\eta-\eta=\eta.$$
Thus, $b_1> a_2$, and also, as $b_2> a_2$ we have,
$$I_1\cap I_2= [a_2, \tilde{b}]\q\h{with}\q  \tilde b > a_2.$$
Next, let  $a_1>a_2$. In this case, $I_1\cap I_2=[a_1,\tilde b]$, where $\tilde b:=\min\{b_2, b_1\}$. Note,  again by (\ref{A-interval}) and (\ref{A-interval-1}), that
$$b_2-a_1 = (b_2-a_2) -(a_1-a_2) \geq 2\eta-\eta=\eta.$$
Thus, $b_2> a_1$, and also, as $b_1> a_1$ we have,
$$I_1\cap I_2= [a_1, \tilde{b}]\q\h{with}\q  \tilde b > a_1.$$
Hence, combining both the cases, we have the proof.
\epf

\brem \label{rem-interval-1}
Let $s_1$ and $s_0$ in $[0,1]$ be such that $g_0=g(\g(s_0))$ and $g_1=g(\g(s_1))$.
Let us recall that $I:= [g_0, g_1]$ and $I_\e:=[g_0^\e, g_1^\e]$. Since $g$ and $g^{\e}$ are in $C^1(\Gamma)$, we have $g\circ\g$ and $g^{\e}\circ\g$ are in $C^1([0,1])$. Also,
$$\|g\circ\g - g^{\e}\circ \g\|_{L^{\infty}([0,1])} \leq \|g-g^{\e} \|_{W^{1,\infty}(\Gamma)} \leq \e.$$
Thus, by Lemma \ref{lem-interval-1}, we have
$$|g_0-g^{\e}_0| < \e \q \h{and} \q |g_1-g^{\e}_1| < \e.$$
Hence, taking $\e < {(g_1-g_0)}/{4}$, we have
\beqarray
(g^{\e}_1- g^{\e}_0) &\geq & |g^{\e}(\g(s_0))-g^{\e}(\g(s_1))|\\
&\geq &   |g_1-g_0|-|g(\g(s_0))-g^{\e}(\g(s_0))|-|g(\g(s_1))-g^{\e}(\g(s_1))|\\
&> & 4{\e}- 2\|g-g^{\e}\|_{W^{1,\infty}(\Gamma)} \\
&> & 4{\e} -2{\e} =2{\e},
\eeqarray
and thus, $2\e < \min\{(g_1-g_0), (g^{\e}_1-g^{\e}_0)\}$. Hence by Lemma \ref{lem-interval-1}, $I\cap I_\e$ is a closed and bounded non-degenerate interval.
Let us denote this interval by $\tilde I_\e$. Thus,
\beq\label{intersection-1}   \tilde I_\e = I\cap I_\e = [\tilde g_0^\e, \tilde g_1^\e]\eeq
for some $\tilde g_0^\e, \tilde g_1^\e \in \R$ with $\tilde g_0^\e <\tilde g_1^\e.$
Also, by Lemma \ref{lem-interval-1} we have,
$$|g_0-\tilde g^\e_0| \leq |g_0-g^{\e}_0| < \e\q\h{and}\q |g_1-\tilde{g}^{\e}_1| \leq |g_1-g^{\e}_1| < \e.\eqno\lozenge$$
\end{remark}

Next, we shall make use of the following lemma whose proof is given in the appendix.

\bl \label{lem-imbedding}
There exists a constant $C>0$ such that for any closed interval $J$,
$$ \|y\|_{L^{\infty}(J)} \leq C_J \|y\|_{H^1(J)},$$
where  $C_J := C \max\{ 3, (2|J|+1)\}$.
In particular, for any interval $J_0$ such that $J_0 \subseteq J$,
\beq \label{ineq-imbedding} \|y\|_{L^{\infty}(J_0)} \leq C_J \|y\|_{H^1(J_0)}.\eeq
\el

If $y\in W^{1,\infty}(J_1)$ then using (\ref{ineq-imbedding}) we obtain
$$ \|y\|^2_{L^{\infty}(J_0)}  \leq  (C_{J_1})^2 \left[ \int_{J_0} y^2 + \int_{J_0} (y')^2 \right]
                           \leq  (C_{J_1})^2 |J_0| \left[ \|y\|^2_{L^{\infty}(J_0)} + \|y'\|^2_{L^{\infty}
                             (J_0)} \right].$$
Thus
\beq \label{ineq-imbedding-1}
\|y\|_{L^2(J_0)} \leq  \sqrt{|J_0|} \|y\|_{L^{\infty}([a,c])} \leq |J_0| \sqrt{2} C_{J_1} \|y\|_{W^{1,\infty}(J_0)},\eeq
and additionally if $y'' \in L^{\infty}(J_1)$, then                      
\beqarray
\|y\|^2_{L^2(J_0)} &\leq & {|J_0|}^2 (C_{J_1})^2 (\|y\|^2_{L^{\infty}(J_0)} + \|y'\|^2_{L^{\infty}(J_0)})\\
                         &\leq & {|J_0|}^3 (C_{J_1})^4 \left[ \|y\|^2_{L^{\infty}(J_0)} + \|y'\|^2_{L^{\infty}
                             (J_0)} + \|y'\|^2_{L^{\infty}(J_0)} + \|y''\|^2_{L^{\infty}
                             (J_0)} \right]\\
                         &\leq & 4{|J_0|}^3 (C_{J_1})^4 \|y\|^2_{W^{2,\infty}(J_0)}
\eeqarray
which implies
\beq \label{ineq-imbedding-2}
\|y\|_{L^2(J_0)} \leq 2{|J_0|}^{3/2} (C_{J_1})^2 \|y\|_{W^{2,\infty}(J_0)}.\eeq

\bl \label{lem-interval-2}
Let $J_1$ and $J_2$ be closed intervals such that $J_2\subseteq J_1$ and let  $C_{J_1}$ be  as in Lemma \ref{lem-imbedding}.
Let $y \in H^2(J_1)$, then  we have the following.
\begin{itemize}
\i[(i)]  $\|y\|_{L^2(J_1\setminus J_2)} \leq  \sqrt{2}C_{J_1}\|y\|_{W^{1,\infty}(J_1)} |J_1 \setminus J_2|. $
\i[(ii)] If $y'' \in L^{\infty}(J_1)$ then
$$ \|y\|_{L^2(J_1\setminus J_2)} \leq  2(C_{J_1})^2\|y\|_{W^{2,\infty}(J_1)} {|J_1 \setminus J_2|}^{3/2}. $$
\end{itemize}
\el

\bpf 
Let $J_1=[a, b]$ and $J_2=[c, d]$ for some $a\leq b$ and $c\leq d$.
If $J_1=J_2$ then $J_1\setminus J_2 =\emptyset$, and in that case the result holds trivially.  So let us consider the cases when either $a<c$ or $d< b$, or both holds.
Without loss of generality let us assume that $a<c$ and $d<b$ .
Let $y\in H^2(J_1)$. Then by (\ref{ineq-imbedding}) $y$ and $y'$ are in $L^{\infty}(J_1)$.
Thus taking $J_0=[a,c]$ in (\ref{ineq-imbedding-1}) we have
 $$ \|y\|_{L^2([a,c])} \leq {(c-a)} \sqrt{2} C_{J_1}\|y\|_{W^{1,\infty}([a,c])} \leq {(c-a)} \sqrt{2}C_{J_1}\|y\|_{W^{1,\infty}(J_1)} $$
and taking $J_0=[d,b]$ in (\ref{ineq-imbedding-1}) we have
$$\|y\|_{L^2([d,b])} \leq (b-d) \sqrt{2} C_{J_1}\|y\|_{W^{1,\infty}([d,b])} \leq {(b-d)} \sqrt{2}C_{J_1}\|y\|_{W^{1,\infty}(J_1)}.$$
Hence we have (i).
Next, additionally if, $y''\in L^{\infty}(J_1)$, having $J_0=[a,c]$ in (\ref{ineq-imbedding-2}) we obtain
$$\|y\|_{L^2([a,c])} \leq 2{(c-a)}^{3/2} (C_{J_1})^2 \|y\|_{W^{2,\infty}([a,c])} \leq 2{(c-a)}^{3/2} (C_{J_1})^2 \|y\|_{W^{2,\infty}(J_1)}$$
and
having $J_0=[d,b]$ in (\ref{ineq-imbedding-2}) we obtain
$$\|y\|_{L^2(d,b)} \leq 2{(b-d)}^{3/2} (C_{J_1})^2 \|y\|_{W^{2,\infty}([d,b])} \leq 2{(b-d)}^{3/2} (C_{J_1})^2 \|y\|_{W^{2,\infty}(J_1)}.$$
Hence we have (ii).
\epf

\bl \label{lem-interval-3}
Let $\phi_1, \, \phi_2, \, I_1,\, I_2$ and $\eta$ be as in Lemma \ref{lem-interval-1} satisfying all the assumptions there. Then, for any interval $I_3 \subset I_1 \cap I_2$ and $y\in C^1(I_1)$
\beq \label{ineq-interval-3-1}\int_{I_3} |y(\phi_1(\xi))- y(\phi_2(\xi))|^2 d\xi \leq {\|y'\|}^2_{L^{\infty}(I_1)} {\|\phi_1 - \phi_2 \|}^2_{L^2([\xi_1, \xi_2])}.\eeq
Assume, further, that  $\phi_1, \phi_2\in C^1([\xi_1, \xi_2])$ satisfying $|{\phi_1}'(\xi)| \geq C_{\phi_1}$ and
$|{\phi_2}'(\xi)| \geq C_{\phi_2}$ for some constants $C_{\phi_1}, C_{\phi_2}>0 $. Then, for $y\in H^2(I_1)$
\beq \label{ineq-interval-3-2}{\|y\circ \phi_1- \widetilde{y \circ \phi_2}\|}^2_{L^2([\xi_1, \xi_2])} \leq C_I \|y\|_{H^2(I_1)} \Big({\|\phi_1 - \phi_2 \|}^2_{L^2([\xi_1, \xi_2])}
+ \frac{2\sqrt{2}}{\sqrt{C_{\phi_1}}} {\eta}\Big),\eeq
where
$$\widetilde{y\circ \phi_2}(\xi):= \left\{\begin{array}{ll}
                              (y\circ \phi_2)(\xi) &\h{if}\q  \xi\in [\tilde{\xi_1}, \tilde{\xi_2}],\\
                              0 &\h{if}\q  \xi\in [{\xi_1}, {\xi_2}]\setminus [\tilde{\xi_1}, \tilde{\xi_2}],\end{array}\right.$$
with $[\tilde{\xi_1}, \tilde{\xi_2}]=(\phi_2)^{-1}(I_1 \cap I_2 )$ and $C_I$ is as in Lemma \ref{lem-imbedding}.
\el

\bpf
By Lemma \ref{lem-interval-1} we have $I_1 \cap I_2$ to be a closed non-degenerate interval. Let $I_3$ be an interval in $I_1 \cap I_2$. Then for $y\in C^1(I_1)$ using fundamental theorem of calculus and H\"{o}lder's inequality we have
\beqarray
\int_{I_3} |y(\phi_1(\xi))- y(\phi_2(\xi))|^2 d\xi
&=& \int_{I_3} {\left[ \int_{\phi_1(\xi))} ^{\phi_2(\xi))} y'(\theta) d\theta \right]}^2 d\xi \\
& \leq & \int_{I_3} {\|y'\|}^2_{L^{\infty}(I_1)}|\phi_1(\xi)) - \phi_2(\xi))|^2 d\xi \\
& \leq &  {\|y'\|}^2_{L^{\infty}(I_1)} {\|\phi_1 - \phi_2 \|}^2_{L^2([\xi_1, \xi_2])}.
\eeqarray
Hence we have (\ref{ineq-interval-3-1}). Next, let $I_1 \cap I_2$ be equal to $[\tilde{a_2}, \tilde{b_2}]$ for some $\tilde{a_2}$ and $\tilde{b_2}$ in $\R$, with $\tilde{a_2}< \tilde{b_2}$.
Since $\phi_2\in C^1([\xi_1, \xi_2])$ and $|{\phi_2}'(\xi)| \geq C_{\phi_2}>0$, $\phi_2$ is invertible from its image and the inverse is continuous. Thus $(\phi_2)^{-1}(I_1 \cap I_2 )= [\tilde{\xi_1}, \tilde{\xi_2}]$ for some $[\tilde{\xi_1},  \tilde{\xi_2}]\subseteq [\xi_1, \xi_2]$.
Also, by the properties of $\phi_1$, we have, $\phi_1([\tilde{\xi_1}, \tilde{\xi_2}])=[\tilde{a_1}, \tilde{b_1}]$ for some $\tilde{a_1}$ and $\tilde{b_1}$ in $I_1$, with $\tilde{a_1} < \tilde{b_1}$.
Thus using Lemma \ref{lem-interval-1} with $[\tilde{a_1}, \tilde{b_1}]$ and $[\tilde{a_2}, \tilde{b_2}]$ in place of $I_1$ and $I_2$ respectively, we have $|\tilde{a_1}- \tilde{a_2}| \leq \eta$ and $|\tilde{b_1}- \tilde{b_2}| \leq \eta $. Hence, using Lemma \ref{lem-interval-1} and definition of $\tilde{a_2}$ and $\tilde{b_2}$  we have,
\beq \label{eq-interval-3} |a_1-\tilde{a_1}| \leq |a_1-\tilde{a_2}|+|\tilde{a_2}- \tilde{a_1}| \leq |a_1-a_2|+
|\tilde{a_2}- \tilde{a_1}| \leq 2\eta, \eeq
\beq \label{eq-interval-4} |b_1-\tilde{b_1}| \leq |b_1-\tilde{b_2}|+|\tilde{b_2}- \tilde{b_1}| \leq
|b_1-b_2|+|\tilde{b_2}- \tilde{b_1}| \leq 2\eta. \eeq
Thus by definition of $\widetilde{y_\circ \phi_2}$, we have
\beq \label{int-3-1}
{\|y\circ \phi_1- \widetilde{y_\circ \phi_2}\|}^2_{L^2([\xi_1, \xi_2])} =
 \int_{\tilde{\xi_1}} ^{\tilde{\xi_2}} |y(\phi_1(\xi))- y(\phi_2(\xi))|^2 d\xi + \int_{[\xi_1, \xi_2]\setminus [\tilde{\xi_1}, \tilde{\xi_2}]} |y(\phi_1(\xi))|^2 d\xi. \eeq
For any $\xi \in [\xi_1, \xi_2]$, $|\phi'_1(\xi)| \geq C_{\phi_1}$ hold. Thus, by Lemma \ref{lem-meas},
 $$\int_{[\xi_1, \xi_2]\setminus [\tilde{\xi_1}, \tilde{\xi_2}]} |y(\phi_1(\xi))|^2 d\xi
  \leq  \frac{1}{C_{\phi_1}} \int_{I_1\setminus \phi_1([\tilde{\xi_1}, \tilde{\xi_2}])} |y(z)|^2 dz.$$
Hence, as (\ref{eq-interval-3}) and (\ref{eq-interval-4}) hold and (\ref{A-interval}) is assumed, taking $J_1=I_1$ and $J_2=\phi_1([\xi_1, \xi_2])=[\tilde{a_1}, \tilde{b_1}]$ in Lemma \ref{lem-interval-2}-(i)
we obtain
$$
  \int_{[\xi_1, \xi_2]\setminus [\tilde{\xi_1}, \tilde{\xi_2}]} |y(\phi_1(\xi))|^2 d\xi
  \leq  \frac{1}{C_{\phi_1}} {\|y\|}^2_{L^2(I_1\setminus \phi_1([\tilde{\xi_1}, \tilde{\xi_2}]))}\\
 \leq  \frac{2(C_{I_1})^2}{C_{\phi_1}} {\|y\|}^2_{H^2(I_1)} {4\eta}^{2}.$$
Thus using (\ref{int-3-1}), the fact that $H^2(I_1)$ is continuously imbedded in $C^1(I_1)$ and having $I_3=[\tilde{\xi_1}, \tilde{\xi_2}]$ in (\ref{ineq-interval-3-1}) we obtain
$$ {\|y\circ \phi_1- \widetilde{y_2\circ \phi_2}\|}_{L^2([\xi_1, \xi_2])} \leq
(\|y'\|_{L^{\infty}(J_1)}\|\phi_1 - \phi_2 \|_{L^2([\xi_1, \xi_2])} + \|y\|_{H^2(J_1)}\frac{2\sqrt{2}C_{I_1}}{\sqrt{C_{\phi_1}}} {\eta}) . $$
Hence, using (\ref{ineq-imbedding}) we have (\ref{ineq-interval-3-2}).
\epf

Let us recall that $I=g(\gamma([0,1]))=[g_0,g_1]$, $I_{\e}=g^{\e}(\gamma([0,1]))=[g^{\e}_0, g^{\e}_1]$ and for $\e< {(g_1-g_0)}/{4}$, let  $\tilde I_\e = I\cap I_\e = [\tilde g_0^\e, \tilde g_1^\e]$ as in  (\ref{intersection-1}).
By (\ref{pert-1}) we have $\|g-g^{\e}\|_{W^{1,\infty}(\Gamma)} \leq \e$ and thus
\beq \label{pert-5}\|g\circ\gamma - g^{\e}\circ\gamma\|_{L^{\infty}([0,1])}=\sup_{s\in [0,1]}|g\circ\gamma(s)- g^{\e}\circ \gamma(s)| \leq \|g-g^{\e}\|_{W^{1,\infty}(\Gamma)} \leq \e.\eeq
Now, additionally let $\e \leq {C_g}/{2}$. Then, by (\ref{ineq-pert}) and (\ref{A-1}) $g^{\e}$ and $\g$ are  bijective, and so $(g^{\e}\circ\gamma)^{-1}$ is continuous. Thus $(g^{\e}\circ\gamma)^{-1}(\tilde{I_{\e}})$ is a closed non-degenerate interval. In other words
\beq \label{defn-interval} \tilde{I_{\e}}=[\tilde{g^{\e}_0},\tilde{g^{\e}_1}]= g^{\e}(\gamma([t^{\e}_0,t^{\e}_1])\eeq
for some $t^{\e}_0$ and  $t^{\e}_1$ in $[0,1]$ with $t^{\e}_0 < t^{\e}_1$.

Now, for $\e \leq \min\{{(g_1-g_0)}/{4},\,\, {C_g}/{2}\}$, let $T_3^{\e}:L^2(\tilde{I_{\e}}) \to L^2([0,1])$ be defined by
\beq \label{op-defn-1} T_3^{\epsilon}(\zeta)(s)=   \left\{
\begin{array}{ll}
       \zeta( g^\e( \gamma(s)))  & s\in [t^{\e}_0, t^{\e}_1]  \\
       0 & g^{\e}(\gamma(s)) \in  [0,1]\setminus [t^{\e}_0, t^{\e}_1].
\end{array}
\right. \eeq

Now, we prove some properties of $T_3^{\e}$.

\bt \label{lem-pert}
Let $T_3^{\e}$ be as defined in (\ref{op-defn-1}). Then, for $\zeta \in \mathcal{W}$,
$$\|T_3\zeta -T_3^{\e}\zeta_{|_{\tilde{I}_\e}}\|_{L^2([0,1])} \leq (C_{g,\gamma, I} \|\zeta\|_{H^2(I)})\e,$$
where $C_{g,\gamma,I}=C_I\Big(1+\frac{2\sqrt{2}}{\sqrt{C_gC_{\gamma}}} \Big)$ with $C_I$ as in (\ref{ineq-imbedding}).
\et

\bpf
Let $\zeta \in \mathcal{W}$.
For any $s\in [0,1]$, by (\ref{A-2}) and (\ref{A-1}), we have
$$|{(g\circ \gamma)}'(s)| \geq C_gC_{\gamma}.$$
By (\ref{pert-5}) and (\ref{defn-interval}), we have  $\|g\circ\gamma -g^{\e}\circ\gamma \| \leq \e$ and $\tilde{I_{\e}}= I \cap I_{\e}= (g^{\e}\circ \gamma)([t^{\e}_0, t^{\e}_1])$, respectively. Now $\zeta \in \mathcal{W} \subset H^2(I)$. Then, by definition of $T_3$ and $T_3^{\e}$, we have
$$T_3(\zeta)=\zeta\circ g\circ \gamma \in L^2([0,1])\q\h{and}\q {(\zeta \circ g^{\e} \circ \gamma)} _{|_{[t^{\e}_0, t^{\e}_1]}}= ({T_3^{\e}(\zeta_{|_{\tilde{I}_\e}})})_{|_{[t^{\e}_0, t^{\e}_1]}} \in L^2([t^{\e}_0, t^{\e}_1]).$$
Hence, taking $\phi_1$ as $g\circ\gamma$ and $\phi_2$ as $g^{\e}\circ\gamma$ in Lemma \ref{lem-interval-3}, we have
$$\|T_3\zeta-T_3^{\e}\zeta_{\e}\|_{L^2([0,1])} =  \|\zeta\circ g\circ \gamma - {(\zeta \circ g^{\e} \circ \gamma)} _{|_{[t^{\e}_0, t^{\e}_1]}} \|_{L^2([0,1])}
\leq  C_I\Big(1+\frac{2\sqrt{2}}{\sqrt{C_gC_{\gamma}}} \Big) \|\zeta\|_{H^2(I)} \e.$$
This completes the  proof.
\epf

\bt \label{Th-T_3_eps_bdd_below}
The map $T_3^{\e}:L^2(\tilde{I_{\e}}) \to L^2([0,1])$, defined as in (\ref{op-defn-1}),  is bounded linear and bounded below. In fact, for every $\zeta\in L^2(\tilde{I_{\e}})$,
\beq \label{T_3-pert-bdd-below-2} \sqrt{\frac{C_g C_\gamma}{2}} \|T^{\e}_3(\zeta)\|_{L^2([0,1])} \leq \|\zeta\|_{L^2(\tilde{I_{\e}})} \leq \sqrt{2C'_g C'_\gamma} \|T^{\e}_3(\zeta)\|_{L^2([0,1])}, \eeq
where $C_\g, C_\g'$ and  $C_g, C_g'$ are as in (\ref{A-1}) and (\ref{A-2}), respectively.\et
\bpf
Clearly, $T^{\e}_3$ is a linear map. Since (\ref{ineq-pert}) and (\ref{A-1}) hold, using Lemma \ref{lem-meas}, and (\ref{op-defn-1}) we obtain
\beqarray
\|T^{\e}_3(\zeta)\|_{L^2([0,1])} & = & \int_0^1 {|T^{\e}_3(\zeta)(s)|}^2 ds  =  \int_{t^{\e}_0} ^{t^{\e}_1} {|\zeta(g^{\e}(\gamma(s)))|}^2 ds\\
&\leq & \frac{2}{C_gC_{\gamma}}\int_{\tilde{g^{\e}_0}}^{\tilde{g^{\e}_1}} |\zeta(z)|^2 dz
= \frac{2}{C_gC_\gamma} \|\zeta\|_{L^2([{\tilde{g^{\e}_0}}, {\tilde{g^{\e}_1}}])},
\eeqarray
\beqarray
\|T^{\e}_3(\zeta)\|_{L^2([0,1])} & = & \int_0^1 {|T^{\e}_3(\zeta)(s)|}^2 ds
 =  \int_{t^{\e}_0} ^{t^{\e}_1} {|\zeta(g^{\e}(\gamma(s)))|}^2 ds\\
&\geq & \frac{1}{2C'_gC'_{\gamma}}\int_{\tilde{g^{\e}_0}}^{\tilde{g^{\e}_1}} |\zeta(z)|^2 dz
= \frac{1}{2C'_gC'_{\gamma}} \|\zeta\|_{L^2([{\tilde{g^{\e}_0}}, {\tilde{g^{\e}_1}}])}.
\eeqarray
Hence we have the proof.
\epf

Now, by Theorem \ref{Th-T_3_eps_bdd_below} we know that $T^{\e}_3$ is a bounded linear map which is bounded below. Thus using Lemma \ref{lem-bdd-bel}, the operator
$$(T_3^\e )^\dagger:= {({(T_3^{\e})}^*T_3^{\e})}^{-1} (T_3^{\e})^*$$ is a bounded linear operator and is the generalized inverse of $T_3^\e$.  The following theorem, which also follows from Lemma \ref{lem-bdd-bel}, shows that the family
$$\left \{(T_3^\e )^\dagger: 0 <\e \leq \min \{\frac{C_g}{2},\frac{g_1-g_0}{4}\}\, \right\}$$
is in fact uniformly bounded.

\bt For every  $\zeta \in L^2([0,1])$,
 \beq \label{T_3-inv}\|(T_3^\e )^\dagger\zeta\|_{L^2(\tilde{I_\e})} \leq \sqrt{2 C'_g C'_\gamma} \|\zeta\|_{L^2([0,1])},\eeq
 where $C'_g$ and $C'_{\gamma}$ are as in (\ref{A-1}) and (\ref{A-2}).
\et

In order to obtain an approximate solution of (\ref{ill-1}) under the nosy data $(j^\d, g^\e)$ satisfying (\ref{pert-1}) and (\ref{pert-2}), we adopt the following operator procedure:
First we  consider the following operator equation
\beq \label{first-op-eq} {(T^{\e}_3)}^*(T^{\e}_3)\zeta= {(T^{\e}_3)}^*f^{j^{\delta}}.\eeq
Let $\tilde{\zeta}_{\e,\d}\in L^2(\tilde{I_{\e}})$ be the unique solution of (\ref{first-op-eq}), that is,
$\tilde{\zeta}_{\e,\d}:= (T_3^\e )^\dagger f^{j^{\delta}}$.
Then, we see that 
$$\zeta_{\e,\d} = \left\{\begin{array}{ll}
\tilde{\zeta}_{\e,\d}& \h{ on }\, \tilde{I_{\e}},\\
0& \h{ on }\, I\setminus \tilde{I_{\e}},\end{array}\right.$$
belongs to $L^2(I)$. Next, we consider  the operator equation 
\beq \label{3.3.3} {(T^\a _2)}^* (T^\a _2) (w)= {(T^\a _2)}^*\zeta_{\e, \d}. \eeq
Let $b_{\a, \e, \d}$ be the unique solution of equation (\ref{3.3.3}).
Thus by solving the operator equations (\ref{first-op-eq}) and (\ref{3.3.3}) we obtain $b_{\a,\e,\d}$.
Since  $b_{\a,\e,\d} \in \mathcal{W} \subset R(T_1)$, $a_{\a,\e,\d}:=b'_{\a,\e,\d}$ is  the solution of the equation
$$T_1(a)=b_{\a,\e,\d}.$$
We show that $a_{\a,\e,\d}$ is a candidate for an approximate solution to Problem (P).

\bl \label{lem-err}
Under the assumptions in Assumption \ref{assum-1} on $(j, g)$,
let $a_0 \in H^1(I)$ be the solution of $T(a)=f^j$. Assume further that $a_0(g_1)=0$.
For $\zeta \in L^2(I)$,  let $b_{\a,\zeta} \in H^2(I)$ be such that 
$${(T^\a _2)}^* (T^\a _2) (b_{\a,\zeta})= {(T^\a _2)}^*\zeta,$$
 and let $a_{\a,\zeta}=b'_{\a,\zeta}$. Then
\begin{eqnarray}
\label{lem-1err-est-h1} \|a_0-a_{\a,\zeta} \|_{H^1(I)} 
&\leq&  C_{\a} + \frac{\|\zeta-b_0\|_{L^2(I)}}{\a},\\
\label{lem-1err-est-l2} \|a_0- a_{\a,\zeta} \|_{L^2(I)} 
&\leq & \sqrt{\a}\|a'_0\|_{L^2(I)}
+ \frac{\|\zeta-b_0\|_{L^2(I)}}{ \sqrt{\a}},
\end{eqnarray}
where $C_{\a}>0$ is such that $C_{\a} \to 0$ as $\a\to 0$.
In addition, if $a_0 \in H^3(I)$, then
\begin{eqnarray}
\label{lem-2err-est-h1} \|a_0-a_{\a,\zeta} \|_{H^1(I)} 
&\leq& (1+ C_L)\a \|a'_0\|_{H^2(I)} + \frac{\|\zeta-b_0\|_{L^2(I)}}{\a},\\
\label{lem-2err-est-l2} \|a_0- a_{\a,\zeta} \|_{L^2(I)} 
&\leq & (1+ C_L)\a \|a'_0\|_{H^2(I)}
+ \frac{\|\zeta-b_0\|_{L^2(I)}}{\sqrt{\a}}.
\end{eqnarray}
Here $C_L$ is as (\ref{C_L}).
\el
\bpf
Let $b_0=T_1(a_0)$. Then, as $a_0(g_1)=0$, we have $b_0 \in \mathcal{W}$.
Now, by definition of $a_{\a,\zeta}$ and, $H^1(I)$ and $H^2(I)$ norms, for $r\in \{0,1\}$
\beqarray
\|a_0-a_{\a,\zeta}\|_{H^r(I)}&=& \|a_0-({({(T_2^\a)}^*T_2^\a)}^{-1}{(T_2^\a)}^* \zeta)'\|_{H^r(I)} \\
&\leq & \|b_0-{({(T_2^\a)}^*T_2^\a)}^{-1}{(T_2^\a)}^* \zeta\|_{H^{r+1}(I)}\\
&\leq & \|b_0-{({(T_2^\a)}^*T_2^\a)}^{-1}{(T_2^\a)}^*T_2(b_0)\|_{H^{r+1}(I)} \\
&& + \|{({(T_2^\a)}^*T_2^\a)}^{-1}{(T_2^\a)}^* (\zeta-T_2(b_0))\|_{H^{r+1}(I)}
\eeqarray
Hence, for $r \in \{0,1\}$,
\begin{eqnarray} \label{lem-3.4.1}\|a_0-a_{\a,\zeta}\|_{H^r(I)} &\leq & \|b_0-{({(T_2^\a)}^*T_2^\a)}^{-1}{(T_2^\a)}^*T_2(b_0)\|_{H^{r+1}(I)}\\
&&\nonumber  + \|{({(T_2^\a)}^*T_2^\a)}^{-1}{(T_2^\a)}^* (\zeta-T_2(b_0))\|_{H^{r+1}(I)}. \end{eqnarray}
By Theorem \ref{th-reg-conv} we have
\beq \label{lem-3.4.2}
 \|b_0-{({(T_2^\a)}^*T_2^\a)}^{-1}{(T_2^\a)}^*T_2(b_0)\|_{H^2(I)} \to 0 \q \h{as }\a \to 0. \eeq
Also, by Theorem \ref{th-reg}-(iii) we have
\beq \label{lem-3.4.2.2} \|b_0-{({(T_2^\a)}^*T_2^\a)}^{-1}{(T_2^\a)}^*T_2(b_0)\|_{H^1(I)} \leq \|b''_0\|_{L^2(I)} \sqrt{\a}. \eeq
Again, using (\ref{ineq-reg-1}) and (\ref{ineq-reg-2}), we have
\beq \label{lem-3.4.3} \|{({(T_2^\a)}^*T_2^\a)}^{-1}{(T_2^\a)}^* (\zeta-T_2(b_0))\|_{H^2(I)} \leq
  \frac{1}{ \a} \|\zeta - T_2(b_0)\|_{L^2(I)} \eeq
and
\beq \label{lem-3.4.3.2} \|{({(T_2^\a)}^*T_2^\a)}^{-1}{(T_2^\a)}^* (\zeta-T_2(b_0))\|_{H^1(I)} \leq
  \frac{1}{ \sqrt{\a}} \|\zeta - T_2(b_0)\|_{L^2(I)}. \eeq
Thus combining (\ref{lem-3.4.1}), (\ref{lem-3.4.2}) and (\ref{lem-3.4.3}) we have (\ref{lem-1err-est-h1}) with
$$C_\a:= \|b_0-{({(T_2^\a)}^*T_2^\a)}^{-1}{(T_2^\a)}^*T_2(b_0)\|_{H^2(I)},$$ 
 and combining (\ref{lem-3.4.1}), (\ref{lem-3.4.2.2}) and (\ref{lem-3.4.3.2}) we have (\ref{lem-1err-est-l2}).

Next, let $a_0 \in H^3(I)$, $b_0=T_1(a_0) \in \mathcal{W}\cap H^4(I)$.
Then, using theorem \ref{th-reg}-(ii) we have, for $r \in \{0,1\}$,
\beq \label{lem-3.5.2}
 \|b_0-{({(T_2^\a)}^*T_2^\a)}^{-1}{(T_2^\a)}^*T_2(b_0)\|_{H^{r+1}(I)} \leq (1+ C_L) \|b''_0\|_{H^2(I)} \a. \eeq
Thus combining (\ref{lem-3.4.1}), (\ref{lem-3.5.2}) and (\ref{lem-3.4.3}) we have (\ref{lem-2err-est-h1}), and
combining (\ref{lem-3.4.1}), (\ref{lem-3.5.2}) and (\ref{lem-3.4.3.2}) we have (\ref{lem-2err-est-l2}).
\epf

Now, we prove one of the main theorems of this paper.

\bt \label{Th-err}
Let $\e< \min \{{(g_1-g_0)}/{4},\,\, {C_g}/{2}\}$.
Let $a_0$, $g$ and $j$ be as in Lemma \ref{lem-err}.
Let $g^{\e}\in C^1(\Gamma)$, $j^{\d} \in W^{1-1/p,p}(\partial \O)$ with $p>3$,
$\zeta_{\e,\d}$ be the solution of (\ref{first-op-eq}), and $a_{\a,\e,\d}=b'_{\a,\e,\d}$ where $b_{\a,\e,\d}$ is the solution of (\ref{3.3.3}).
Also, let $g^{\e}$ and $j^{\d}$ satisfy (\ref{pert-1}) and (\ref{pert-2}), respectively. Then
%
\begin{eqnarray} \label{1err-est-h1} 
\q\q\q\|a_0-a_{\a,\e,\d} \|_{H^1(I)} &\leq & C_{\a} + \frac{1}{\a}[\sqrt{2C'_g C'_\gamma} ( C_{I,g,\gamma} \|b_0\|_{H^2(I)} \e + \tilde{C_\gamma}\d)+C_I \|b_0\|_{H^2(I)}{\e}  ],\\
\label{1err-est-l2}\q\q\q \|a_0- a_{\a,\e,\d} \|_{L^2(I)} &\leq& \sqrt{\a}\|a'_0\|_{L^2(I)}
+ \frac{1}{ \sqrt{\a}}[\sqrt{2C'_g C'_\gamma} ( C_{I,g,\gamma} \|b_0\|_{H^2(I)} \e + \tilde{C_\gamma}\d)+C_I \|b_0\|_{H^2(I)}{\e}  ], \end{eqnarray}
where $C_{\a}>0$ is such that $C_{\a} \to 0$ as $\a\to 0$.

In addition if $a_0 \in H^3(I)$, then
\beq \label{2err-est-h1} \|a_0-a_{\a,\e,\d} \|_{H^1(I)} \leq (1+ C_L)\|a'_0\|_{H^2(I)}\a + \frac{1}{\a}
[\sqrt{2C'_g C'_\gamma} ( C_{I,g,\gamma} \|b_0\|_{H^2(I)} \e + \tilde{C_\gamma}\d)+C_I \|b_0\|_{H^2(I)}{\e}], \eeq
\beq \label{2err-est-l2} \|a_0- a_{\a,\e,\d} \|_{L^2(I)} \leq (1+ C_L)\|a'_0\|_{H^2(I)}\a
+ \frac{1}{\sqrt{\a}}[\sqrt{2C'_g C'_\gamma} ( C_{I,g,\gamma} \|b_0\|_{H^2(I)} \e + \tilde{C_\gamma}\d)+C_I \|b_0\|_{H^2(I)}{\e}  ]. \eeq
Here, $b_0=T_1(a_0)$, and $C_L$, $\tilde{C_{\gamma}}$, $C_I$, $C_{I,g,\gamma}$, $C'_g$ and $C'_\gamma $ are constants as in (\ref{C_L}), Lemmas \ref{p-trace} and \ref{lem-imbedding}, Theorem \ref{lem-pert}, (\ref{A-1}) and (\ref{A-2}) respectively.
\et
\bpf
Since $a_0(g_1)=0$, we have $b_0 \in \mathcal{W}$.
Now let us note that, by Remark \ref{rem-interval-1}, we have
$|g_0-\tilde{g^{\e}_0}| < \e$ and $|g_1-\tilde{g^{\e}_1}| < \e$.
 Hence, taking $J_1$ and $J_2$ as $I$ and $\tilde{I_{\e}}$ respectively in Lemma \ref{lem-interval-2}, and with our choice of $\e$, by Lemma \ref{lem-interval-2}-(i) we have,
 \beq \label{ineq-err-1}\|b_0\|_{L^2(I\setminus \tilde{I_{\e}})} \leq C_{I}\|b_0\|_{H^2(I)}{\e}. \eeq
Since  $\tilde{\zeta}_{\e,\d}=  {(T^{\e}_3)}^\dagger f^{j^{\delta}}$, $T_3(T_2(b_0))=f^j$, and  ${(T^{\e}_3)}^\dagger T_3^\e $ is identity, we have
\beqarray
\|\tilde \zeta_{\e,\d} - {b_0}|_{\tilde{I_{\e}}}\|_{L^2(\tilde{I_{\e}})} &=& \|\tilde{\zeta}_{\e,\d} - (T_2(b_0))|_{\tilde{I_{\e}}}\|_{L^2(\tilde{I_{\e}})}\\
 &= & \| {(T^{\e}_3)}^\dagger f^{j^{\delta}} -(T_2(b_0))|_{\tilde{I_{\e}}} \|_{L^2(\tilde{I}_\e)} \\
 &\leq & \| (T_3^\e )^\dagger  T_3(T_2(b_0))-(T_2(b_0))|_{\tilde{I_{\e}}}\|_{L^2(\tilde{I}_\e)}  +
 \|(T_3^\e )^\dagger  (f^{j}-f^{j^\d})\|_{L^2(\tilde{I}_\e)}\\
 &\leq & \|(T_3^\e )^\dagger  (T_3(T_2(b_0))-T_3^{\e}((T_2(b_0))|_{T_2(b_0)}))\|_{L^2(\tilde{I}_\e)} + \|(T_3^\e )^\dagger (f^{j}-f^{j^\d})\|_{L^2(\tilde{I}_\e)}.
\eeqarray
Now, by (\ref{T_3-inv}) and Theorem \ref{lem-pert}, we obtain
$$ \|(T_3^\e )^\dagger (f^{j}-f^{j^\d})\|_{L^2(\tilde{I}_\e)}\leq \sqrt{2C'_g C'_\gamma} \tilde{C_\gamma}\d,$$
$$
\|(T_3^\e )^\dagger (T_3(T_2(b_0))-T_3^{\e}((T_2(b_0))|_{T_2(b_0)}))\|_{L^2(\tilde{I}_\e)}
\leq  \sqrt{2C'_g C'_\gamma} C_{I,g,\gamma} \|b_0\|_{H^2(I)} \e.$$
Therefore,
 \beq \label{ineq-err-2} \|\tilde{\zeta}_{\e,\d} - {b_0}|_{\tilde{I_{\e}}}\|_{L^2(\tilde{I_{\e}})} \leq \sqrt{2C'_g C'_\gamma} [ (C_{I,g,\gamma}) \|b_0\|_{H^2(I)} \e + \tilde{C_\gamma}\d]. \eeq
Now by definition of $\zeta_{\e,\d}$ we have
$$\|\zeta_{\e,\d}-b_0\|_{L^2(I)} \leq \|\tilde{\zeta}_{\e,\d}-{b_0}|_{\tilde{I_{\e}}}\|_{L^2(\tilde{I_\e})} + \|b_0\|_{L^2(I\setminus \tilde{I_\e})}.$$ Hence, by (\ref{ineq-err-1}) and (\ref{ineq-err-2}) we have,
$$\|\zeta_{\e,\d}-b_0\|_{L^2(I)} \leq \sqrt{2C'_g C'_\gamma}   [ (C_{I,g,\gamma}) \|b_0\|_{H^2(I)} \e + \tilde{C_\gamma}\d] +C_I \|b_0\|_{H^2(I)}{\e}.$$
Now by definition, $b_{\a,\e,\d}$ is the unique solution of equation (\ref{3.3.3}).
Thus, with $\zeta_{\e,\d}$ in place of $\zeta$ in Lemma \ref{lem-err}, we have the proof.
\epf

\brem \label{rem-err-est-1}
Let $a_0$ and $a_{\a,\e,\d}$ be as defined in Theorem \ref{Th-err}.  Then
(\ref{1err-est-h1}) and (\ref{1err-est-l2}) take the forms
\beqarray
\|a_0-a_{\a,\e,\d} \|_{H^1(I)} &\leq&  C_{\a} + K_1\, \frac{\e+\d}{\a},\\
\|a_0- a_{\a,\e,\d} \|_{L^2(I)} &\leq & \sqrt{\a}\|a'_0\|_{L^2(I)} + K_2\, \frac{\e+\d}{\sqrt{\a}},\eeqarray
respectively, where $C_{\a}>0$ is such that $C_{\a} \to 0$ as $\a\to 0$, and if, in addition, $a_0 \in H^3(I)$, then
(\ref{2err-est-h1}) and (\ref{2err-est-l2}) take the forms
\beqarray
\|a_0-a_{\a,\e,\d} \|_{H^1(I)} &\leq &   (1+ C_L)\|a'_0\|_{H^2(I)}\a + K_3\, \frac{\e+\d}{\a},\\
\|a_0- a_{\a,\e,\d} \|_{L^2(I)} &\leq & (1+ C_L)\|a'_0\|_{H^2(I)}\a + K_4\,  \frac{\e+\d}{\sqrt{\a}},\eeqarray
respectively, where $K_1, K_2, K_3, K_4$ are positive constants independent of $\a, \e, \d$ and $C_L\geq \|id-L\|$, where $L$ is the bounded operator as in Proposition \ref{prop-projection}. 
.

Then, choosing $\a=\sqrt{\d}$ and $\e= \d$ in (\ref{1err-est-h1}) we have
$$\|a_0- a_{\a,\e,\d}\|_{H^1(I)}=o(1).$$
Thus using the new regularization method we obtain a result better than the order $O(1)$ in \cite{Egger} obtained using Tikhonov regularization.
On choosing  $\a=\d=\e$ in (\ref{1err-est-l2}) we have
$$\|a_0- a_{\a,\e,\d}\|_{L^2(I)}= O(\sqrt{\d}),$$
which is same as the estimate obtained in \cite{Egger}. Next, under the source condition $a_0 \in H^3(I)$ and for $\a=\sqrt{\d}$ and $\e=\d$, (\ref{2err-est-h1}) gives the order as
$$\|a_0- a_{\a,\e,\d}\|_{H^1(I)}=O(\sqrt{\d}).$$
This estimate is similar to a result obtained in \cite{Kugler} with source condition $a_0\in H^4(I)$ and trace of $a_0$ being Lipschitz which is stronger than the source condition needed in our result, whereas  under the same source condition $a_0 \in H^3(I)$, the choice of $\a={\d}^{2/3}$ and $\e=\d$  in  (\ref{2err-est-l2}) gives the rate as
$$\|a_0- a_{\a,\e,\d}\|_{L^2(I)}= O({\d}^{2/3}).$$
This is better than the rate  $O({\d}^{3/5})$ mentioned in \cite{Egger} as the best possible estimate under $L^2(I)$ norm (under realistic boundary condition) using Tikhonov regularization.
\erem

\section{Relaxation of assumption on perturbed data}\label{sec-4}
In the previous section we have carried out our analysis assuming that the perturbed data $g^{\e}$ is in $C^1(\Gamma)$, along with (\ref{pert-1}). This assumption can turn out to be too strong for implementation in practical problems. Hence, here we consider a weaker and practically relevant assumption on our perturbed data $g^{\e} $, namely $g^{\e} \in L^2(\Gamma)$ with
\beq \label{new-pert-1} \|g-g^{\e}\|_{L^2(\Gamma)} \leq \e. \eeq

What we essentially used in our analysis  in Section \ref{sec-3} to derive the error estimates is that $g^\e\circ \g$ is close to $g\circ \g$ in appropriate norms.   Here, we consider $\tilde g_\g^\e := \Pi_h(g^\e\circ \g)$ in place of $g^\e\circ \g$, where $\Pi_h: L^2([0,1]) \to L^2([0,1])$ is the orthogonal projection onto a subspace of  $W^{1,\infty}([0, 1])$, and we show that
$\tilde g_\g^\e$ is close to $g\circ \g$ in  appropriate norms, and then obtain associated error estimates.
For this purpose, we shall also assume more regularity on $g\circ \g$, namely, $g\circ \g\in H^4([0, 1])$.

Let  $\Pi_h: L^2([0,1]) \to L^2([0,1])$ be the orthogonal projection onto the space $L_h$ which is the space of all continuous real valued piecewise linear functions $w$ on $[0,1]$ defined on a uniform partition
$0=t_0<t_1< \cdots t_N=1$ of mesh size $h$, that is, $t_i := (i-1)h$ for $i=1,\cdots N$ and  $h=1/N$.
Thus, $w\in L_h$ if and only if $w \in C[0,1]$  such that $w_{|_{[t_{i-1}, t_i]}}$ is a polynomial of degree at most $1$.
Let $\mathbb{T}_h := \{[t_{i-1}, t_i] : i=1,\cdots (\frac{1}{h} +1)\} $.

In the following, for $w\in L^2([0,1])$, we use the notation
$\|w\|_{H^m(\tau_h)}$ and $\|w\|_{W^{m,\infty}(\tau_h)}$ whenever $w_{|_{\t_h}}$ belong to  ${H^m(\tau_h)}$ and ${W^{m,\infty}(\tau_h)}$, respectively. As a particular case of inverse inequality stated in Lemma 4.5.3 in \cite{Bren}, for $m \in \{0,1\}$, we have
\beq\label{inv-ineq-1} \|\Pi_h w\|_{W^{m,\infty}(\tau_h)} \leq C_{m}'\frac{1}{h^{(1/2+m)}} \|\Pi_hw\|_{L^2(\tau_h)},\eeq
where $C_m'$ is a positive constant.

\begin{Proposition}
Let  $w\in L^2([0,1])$, $m\in \N \cup \{0\}$  and $\tau_h \in \mathbb{T}_h$. Then the following inequalities hold.
\begin{eqnarray}     \|w\|_{H^m(\tau_h)} &\leq & h^{1/2} C_0 \|w\|_{H^{m+1}(\tau_h)}\q\h{whenever}\q  w_{|_{\t_h}} \in H^{m+1}(\tau_h), \label{eq-proj-1}\\
 \|w\|_{W^{m,\infty}(\tau_h)} &\leq & C^2_0 h^{1/2}\|w\|_{H^{m+2}(\tau_h)} \q\h{whenever}\q w_{|_{\t_h}} \in W^{m,\infty}(\tau_h) \label{eq-proj-2},\\
  \|\Pi_hw \|_{W^{m,\infty}(\tau_h)} &\leq &   C'_m C_0^{(2+2m)} h^{1/2} \|w\|_{H^{(2m+2)}(\tau_h)}\q\h{whenever}\q w_{|_{\t_h}} \in H^{2m+2}(\tau_h),\label{eq-proj-3}
\end{eqnarray}
where $C_0:= 2C_{[0,1]}$ with $C_{[0,1]}$ as in $(\ref{ineq-imbedding})$
and $C_m'$ is as in (\ref{inv-ineq-1}).
\end{Proposition}

\bpf
If  $w^{(j)}_{|_{\t_h}}\in H^1(\tau_h)$ for some $j \in \N \cup \{0\}$, then using (\ref{ineq-imbedding}) and  the fact that $\tau_h$ is of length $h$, we obtain
$$\|w^{(j)}\|_{L^2(\tau_h)}  \leq  h^{1/2} \|w^{(j)}\|_{L^{\infty}(\tau_h)} \leq h^{1/2} C_{I_0} \|w^{(j)}\|_{H^1(\tau_h)},$$
where $I_0:=[0, 1]$.
Hence, we have
$$ \|w\|_{H^m(\tau_h)} = \sum_{j=0}^{m} \|w^{(j)}\|_{L^2(\tau_h)} \leq \sum_{j=0}^{m} h^{1/2} C_{I_0} \|w^{(j)}\|_{H^1(\tau_h)} \leq 2C_{I_0} h^{1/2} \|w\|_{H^{(m+1)}(\tau_h)}. $$
Thus,  taking $C_0=2C_{I_0}$, we  have (\ref{eq-proj-1}).

By repeatedly using (\ref{ineq-imbedding}) and then by (\ref{eq-proj-1}), we obtain
$$\|w\|_{{W^{m,\infty}(\tau_h)}} \leq 2C_{I_0} \|w\|_{H^{m+1}(\tau_h)} \leq 2C_{I_0}C_0 h^{1/2}\|w\|_{H^{m+2}(\tau_h)}. $$
As we have taken $C_0=2C_{I_0}$, we have the proof of (\ref{eq-proj-2}).

Since $\Pi_h$ is an orthogonal projection, from (\ref{inv-ineq-1}) we obtain,
$$\|\Pi_h w\|_{W^{m,\infty}(\tau_h)}\leq \frac{C_m'}{h^{(1/2+m)}}\|\Pi_h w\|_{L^2(\tau_h)} \leq \frac{C_m'}{h^{(1/2+m)}}\|w\|_{L^2(\tau_h)}, $$
and, by repeatedly using (\ref{eq-proj-1}) we have
$$\frac{C_m'}{h^{(1/2+m)}}\|w\|_{L^2(\tau_h)} \leq  C_0^{(2+2m)}\frac{C_m'}{h^{(1/2+m)}}h^{((2m+2)/2)}
\|w\|_{H^{(2m+2)}(\tau_h)} \leq  C_0^{(2+2m)}C_m' h^{1/2}\|w\|_{H^{(2m+2)}(\tau_h)}.$$
Hence we have the proof of (\ref{eq-proj-3}).
\epf

For simplifying the notation, we shall denote
$$g_\g:= g\circ\gamma,\q g_\g^\e\:= g^\e\circ\gamma.$$
By definition,  $\Pi_h(g_\g^{\e})\in W^{1,\infty}([0,1])$.  In order to show that $\Pi_h(g_\g^{\e})$ is close to $g_\g$ with respect to appropriate norms, we assume that
\beq \label{A-3} g_\g\in H^4([0,1]). \eeq

\bt \label{th-new-pert-1}
Let $\tau_h \in \mathbb{T}_h$ and (\ref{A-3}) be satisfied. Then, the following inequalities hold.
\ben
\i[(i)]
$\|\Pi_h g_\g^\e-g_\gamma \|_{L^{\infty}(\tau_h)} \leq \tilde{C_0} h^{3/2} \|g_\gamma\|_{H^{4}(\tau_h)} +\frac{C'_0}{C_\gamma} \frac{\e}{h^{1/2}}, $

\i[(ii)]
$\|\Pi_h g^\e_\gamma-g_\gamma \|_{W^{1,\infty}(\tau_h)} \leq \tilde{C_1} h^{1/2} \|g_\gamma\|_{H^{4}(\tau_h)} +\frac{C'_1}{C_{\gamma}} \frac{\e}{h^{3/2}}. $

\i[(iii)]
$|(\Pi_h g^\e_\gamma )'(s)| \leq C'_gC'_{\gamma} + \tilde{C_1} h^{1/2} \|g_\gamma\|_{H^{4}(\tau_h)} + \frac{C'_1}{C_{\gamma}} \frac{\e}{h^{3/2}}, \q\forall \, s\in \tau_h,$

\i[(iv)]
$|(\Pi_h g^\e_\gamma)'(s)| \geq C_gC_{\gamma}- \tilde{C_1} h^{1/2} \|g_\gamma\|_{H^{4}(\tau_h)} - \frac{C'_1}{C_{\gamma}} \frac{\e}{h^{3/2}}\q\forall\, s\in \tau_h.$
\een
\et
\bpf
Using triangle inequality we have
\begin{eqnarray}
\label{tr-ineq-1}\|\Pi_hg^{\e}_\gamma -g_\gamma \|_{L^{\infty}(\tau_h)} &\leq&
\|\Pi_hg^{\e}_\gamma-\Pi_hg_\gamma \|_{L^{\infty}(\tau_h)} +
\|\Pi_hg_\gamma-g_\gamma \|_{L^{\infty}(\tau_h)},\\
\label{tr-ineq-2}\|\Pi_hg^{\e}_\gamma-g_\gamma \|_{W^{1,\infty}(\tau_h)} &\leq&
\|\Pi_hg^{\e}_\gamma-\Pi_hg_\gamma \|_{W^{1,\infty}(\tau_h)} +
\|\Pi_hg_\gamma -g_\gamma \|_{W^{1,\infty}(\tau_h)}.
\end{eqnarray}
Assumption (\ref{A-1}), Lemma \ref{lem-meas} and (\ref{new-pert-1}) imply
\beq \label{new-pert-mod} C_\g \|g^{\e}_\gamma- g_\gamma \|_{L^2(\tau_h)}\leq \|g^{\e}- g \|_{L^2(\Gamma)}\leq \e \eeq
so that, using (\ref{inv-ineq-1}) and the fact that $\Pi_h$ is an orthogonal projection, we have
\beq \label{tr-ineq-1-1} \|\Pi_hg^{\e}_\gamma -\Pi_hg_\gamma  \|_{L^{\infty}(\tau_h)}\leq C_0'  \frac{1}{h^{1/2}} \|g^{\e}_\gamma- g_\gamma \|_{L^2(\tau_h)}   \leq \frac{C_0'}{C_\g} \frac{\e}{h^{1/2}},\eeq
and
\beq \label{tr-ineq-2-1} \|\Pi_hg^{\e}_\gamma -\Pi_hg_\gamma  \|_{W^{1,\infty}(\tau_h)}\leq C_1'  \frac{1}{h^{3/2}} \|g^{\e}_\gamma- g_\gamma \|_{L^2(\tau_h)}   \leq \frac{C_1'}{C_\g} \frac{\e}{h^{3/2}},\eeq
By (\ref{eq-proj-2}), (\ref{eq-proj-3}) and (\ref{eq-proj-1})
$$ \|\Pi_hg_\gamma -g_\gamma \|_{L^{\infty}(\tau_h)} \leq \|\Pi_hg_\gamma \|_{L^{\infty}(\tau_h)} + \|g_\gamma \|_{L^{\infty}(\tau_h)} \leq 2{(C_0)}^2 h^{1/2} \|g\circ\gamma\|_{H^2(\tau_h)} \leq 2{(C_0)}^4 h^{3/2} \|g\circ\gamma\|_{H^4(\tau_h)}.$$
Hence, using (\ref{tr-ineq-1}) and (\ref{tr-ineq-1-1}), and taking $\tilde{C_0}=2{(C_0)}^4$, we have (i).
By (\ref{eq-proj-2}) and (\ref{eq-proj-3}),
$$\|\Pi_hg_\gamma -g_\gamma \|_{W^{1,\infty}(\tau_h)} \leq \|\Pi_hg_\gamma \|_{W^{1,\infty}(\tau_h)} + \|g_\gamma \|_{W^{\infty}(\tau_h)} \leq {(C_0)}^4 h^{1/2} \|g\circ\gamma\|_{H^4(\tau_h)} + {(C_0)}^2 h^{1/2}\|g\circ\gamma\|_{H^4(\tau_h)}. $$
Hence, using (\ref{tr-ineq-2-1}) and (\ref{tr-ineq-2}), and taking $\tilde{C_1}={(C_0)}^4 + {(C_0)}^2$ we have (ii).

To prove (iii) and (iv), let $s\in [0,1]$. Note that
$$|(g_\gamma)'(s)| - \|\Pi_hg^{\e}_\gamma -g_\gamma \|_{W^{1,\infty}(\tau_h)} \leq
|(\Pi_hg^{\e}_\gamma )'(s)| \leq \|(g_\gamma)'\|_{L^{\infty}(\tau_h)} + \|\Pi_hg^{\e}_\gamma -g_\gamma \|_{W^{1,\infty}(\tau_h)}.$$
Using (\ref{A-1}) and (\ref{A-2}) the above implies
$$C_gC_{\gamma} - \|\Pi_hg^{\e}_\gamma -g_\gamma \|_{W^{1,\infty}(\tau_h)} \leq
|(\Pi_hg^{\e}_\gamma )'(s)| \leq C'_gC'_{\gamma} + \|\Pi_hg^{\e}_\gamma -g_\gamma \|_{W^{1,\infty}(\tau_h)}.$$
Hence using (ii) we have (iii) and (iv).
\epf

From (iii) and (iv) in Theorem \ref{th-new-pert-1} we obtain the following corollary.
\bcor \label{cor-h}
Let $h$ be such that
\beq \label{ineq-h}  \tilde{C_1} h^{2} \|g_\gamma\|_{H^{4}(\tau_h)} +\frac{C'_1}{C_{\gamma}} \e \leq \frac{C_gC_{\gamma}}{2} h^{3/2}.\eeq
Then,
\beq \label{ineq-pert-h} \frac{C_gC_{\gamma}}{2} \leq |(\Pi_hg^{\e}_\gamma )'(s)| \leq 2C'_gC'_{\gamma}. \eeq
\ecor

Since $(g_\gamma)' \neq 0 $, for any $\tau_h \in \mathbb{T}_h$, $g(\gamma(\tau_h))=[g^h_0, g^h_1]$ for some $g^h_0< g^h_1$. Let us denote
\beq \label{interval-new-1} I_h:=[g^h_0, g^h_1],\q   I^h_{\e}:=\Pi_hg^{\e}_{\gamma}([\tau_h]). \eeq

\begin{Proposition} \label{rem-new-interval}
 Let $h$ and $\e$ satisfy (\ref{ineq-h}) and
 \beq \label{A-h-2} \tilde{C_0} h^{2} \|g_\gamma\|_{H^{4}(\tau_h)} +\frac{C'_0}{C_{\gamma}} {\e}  <\frac{h^{3/2}}{2}.\eeq
 Then, for $\tau_h \in \mathbb{T}_h$, $I_h\cap I^h_{\e}$ is a closed interval with non-empty interior, say $\tilde{I^h_{\e}}=[g^{h,\e}_0, g^{h,\e}_1]$ for some $g^{h,\e}_0 < g^{h,\e}_1$, and
\begin{eqnarray}
\label{ineq-new-int} |g^h_0-g^{h,\e}_0| &<& \tilde{C_0} h^{3/2} \|g_\gamma\|_{H^{4}(\tau_h)} +\frac{C'_0}{C_{\gamma}} \frac{\e}{h^{1/2}},\\
\label{ineq-new-int-1}|g^h_1-g^{h,\e}_1| &<& \tilde{C_0} h^{3/2} \|g_\gamma\|_{H^{4}(\tau_h)}  +\frac{C'_0}{C_{\gamma}} \frac{\e}{h^{1/2}}.
\end{eqnarray}
\end{Proposition}

\bpf
Since $h$ satisfies (\ref{ineq-h}), by Corollary \ref{cor-h}, $\Pi_hg^{\e}_{\gamma}$ satisfies (\ref{ineq-pert-h}). Thus $I^h_{\e}$ is a closed non-degenerate interval. So, by Lemma \ref{lem-interval-1}, taking $\phi_1=(g_\gamma)|_{\tau_h}$ and $\phi_2= (\Pi_hg^{\e}_{\gamma})_{|_{\tau_h}}$ we have the following. $I_h \cap I^h_{\e}= [g^{h,\e}_0, g^{h,\e}_1]$ for some
$g^{h,\e}_0 < g^{h,\e}_1$. Also, since (\ref{A-h-2}) is satisfied, we have (\ref{ineq-new-int}) and (\ref{ineq-new-int-1}).
\epf

Let us recall that, in Section \ref{sec-3} we have the perturbed operator $T^{\e}_3$ corresponding to the perturbed data $g^{\e}$. Here, we are working with $\Pi_h(g^{\e}_{\gamma})$. Now, let us define the corresponding operator which shall be used in place of $T^{\e}_3$, so that we can carry out the analysis similar to that of Section \ref{sec-3}. In order to do that, let us first observe the following.

Let $h$ and $\e$ satisfy (\ref{ineq-h}) and (\ref{A-h-2}). Then, by Corollary \ref{cor-h}, $\Pi_hg^{\e}_{\gamma}$ satisfies (\ref{ineq-pert-h}). Thus, $\Pi_hg^{\e}_{\gamma}$ is bijective and, for any $\tau_h \in \mathbb{T}_h$, ${(\Pi_hg^{\e}_{\gamma})}^{-1}$ is continuous on $I^h_{\e}$. Hence, there exists $t^{h,\e}_0$  and $t^{h,\e}_1$ in $\tau_h$ such that
\beq \label{interval-new-2}\tilde{I^h_{\e}} =[g^{h,\e}_0, g^{h,\e}_1]= \Pi_hg^{\e}_{\gamma} ([t^{h,\e}_0, t^{h,\e}_1]).\eeq
For $y\in L^2(\tilde{I_{\e}})$, let
$$S^{h,\epsilon}(y)(s)=   \left\{
\begin{array}{ll}
       y( \Pi_hg^{\e}_{\gamma}(s))  & s\in [t^{h,\e}_0, t^{h,\e}_1]  \\
       0 & s \in  \tau_h \setminus [t^{h,\e}_0, t^{h,\e}_1] \\

\end{array}
\right. $$
and,  let
\beq\label{op-defn-2} (T^{h,\e}_3y)(s) = (S^{h,\e} y) (s) \q \h{for} \q s\in \tau_h,\,\, \tau_h \in \mathbb{T}_h.\eeq
We observe that $T_3^{h,{\e}} : L^2(\tilde{I_{\e}}) \to L^2([0,1])$ is a linear operator. We shall see some of its properties in the next Theorem.

\bt \label{th-new-bdd-below}
Let $h$ and $\e$ satisfy (\ref{ineq-h}) and (\ref{A-h-2}).
Then, the operator $T_3^{h,{\e}}: L^2(\tilde{I_{\e}}) \to L^2([0,1])$ is bounded linear and bounded below. Further, we have the following.
\begin{itemize}
\i[(i)] For $\zeta \in L^2(\tilde{I_{\e}})$,
\beq \label{new-bdd} \|T_3^{h,{\e}}(\zeta)\|_{L^2([0,1])} \leq \sqrt{\frac{2}{C_gC_{\gamma}}} \|\zeta\|_{L^2(\tilde{I_{\e}})},\eeq
 \beq \label{new-bdd-below} \sqrt{2C'_g C'_\gamma}\|T_3^{h,{\e}}(\zeta)\|_{L^2([0,1])} \geq  \|\zeta\|_{L^2(\tilde{I_{\e}})},\eeq
\beq \label{new-T_3-inv}  \|{({(T_3^{h,\e})}^*T_3^{h,\e})}^{-1} (T_3^{h,\e})^*\| \leq \sqrt{2 C'_g C'_\gamma}.\eeq

\i[(ii)] For $\zeta\in \mathcal{W}$,
\beq \label{new-eq-pert} \|T_3^{h,{\e}}(\zeta_{|_{\tilde{I_{\e}}}})-T_3\zeta\|_{L^2([0,1])} \leq D^1_{g,\e,h}\|\zeta' \|_{H^1(I)} +
\|\zeta\|_{H^2(I)} D^2_{g,\e,h}. \eeq

\i[(iii)] If $\zeta \in \mathcal{W}\cap H^3(I)$, then
\beq \label{new-eq-pert-1} \|T_3^{h,{\e}}(\zeta_{|_{\tilde{I_{\e}}}})-T_3\zeta\|_{L^2([0,1])} \leq D^1_{g,\e,h}\|\zeta' \|_{H^1(I)} +
\|\zeta\|_{H^3(I)} D^3_{g,\e,h}, \eeq
\end{itemize}
where
\beqarray
D^1_{g,\e,h} &=& C_I(\tilde{C_0}\|g\circ\gamma\|_{H^4([0,1])}h^2+ \frac{\e}{C_{\gamma}}),\\
D^2_{g,\e,h} &=&  4(C_I)^2\sqrt{\frac{2}{C_gC_{\gamma}}}{(\tilde{C_0} h^{3/2} \|g_\gamma\|_{H^{4}([0,1])} +\frac{C'_0}{C_{\gamma}} \frac{\e}{h^{1/2}})},\\
D^3_{g,\e,h} &=&  8(C_I)^3\sqrt{\frac{2}{C_gC_{\gamma}}} (\tilde{C_0} h^{3/2} \|g_\gamma\|_{H^{4}([0,1])} +\frac{C'_0}{C_{\gamma}} \frac{\e}{h^{1/2}})^{3/2}.\eeqarray
 Here $C_0$, $\tilde{C_0}$, $C'_0$, $C_I$, $C_g$, $C'_g$, $C'_\gamma $, and $C_\gamma$ are constants as defined in (\ref{eq-proj-1}), Theorem \ref{th-new-pert-1}-(ii), (\ref{ineq-imbedding}), (\ref{A-2}) and (\ref{A-1}) respectively.
\et

\bpf
Clearly, $T^{h,\e}_3$ is a linear map.
Now,
\beq \label{th-new-bdd-below-1}{\|T^{h,\e}_3(\zeta)\|}^2_{L^2([0,1])} =  \sum_{\tau_h \in \mathbb{T}_h} \int_{\tau_h} {|T^{h,\e}_3(\zeta)(s)|}^2 ds
 =  \sum_{\tau_h \in \mathbb{T}_h} \int_{t^{h,\e}_0} ^{t^{h,\e}_1} {|y(\Pi_h g^{\e}_\gamma(s))|}^2 ds. \eeq
Since, $h$ satisfies (\ref{ineq-h}), by Corollary \ref{cor-h}, $\Pi_hg^{\e}_{\gamma}$ satisfies (\ref{ineq-pert-h}). Hence, using Lemma \ref{lem-meas}, we have
\begin{eqnarray}
\label{th-new-bdd-below-2}\sum_{\tau_h \in \mathbb{T}_h} \int_{t^{h,\e}_0} ^{t^{h,\e}_1} {|y(\Pi_h(g^{\e}_\gamma)(s)))|}^2 ds
&\leq & \frac{2}{C_gC_{\gamma}} \sum_{\tau_h \in \mathbb{T}_h} \int_{\tilde{I^h_{\e}}} |\zeta(z)|^2 dz
= \frac{2}{C_gC_{\gamma}} \int_{\tilde{I}_{\e}}  |\zeta(z)|^2 dz,\\
\label{th-new-bdd-below-3} \sum_{\tau_h \in \mathbb{T}_h} \int_{t^{h,\e}_0} ^{t^{h,\e}_1} {|\zeta(\Pi_h(g^{\e}_\gamma)(s)))|}^2 ds
&\geq & \frac{1}{2C'_gC'_{\gamma}} \sum_{\tau_h \in \mathbb{T}_h}  \int_{\tilde{I^h_{\e}}} |\zeta(z)|^2 dz
= \frac{1}{2C'_gC'_{\gamma}} \int_{\tilde{I}_{\e}}  |\zeta(z)|^2 dz.
\end{eqnarray}
Hence, combining (\ref{th-new-bdd-below-1}) and (\ref{th-new-bdd-below-2}) we have (\ref{new-bdd}), and combining (\ref{th-new-bdd-below-1}) and (\ref{th-new-bdd-below-3}) we have (\ref{new-bdd-below}).
Hence, $T^{h,\e}_3$ is bounded linear and bounded below.
Since, $T^{h,\e}_3$ satisfies (\ref{new-bdd}) and (\ref{new-bdd-below}), from Lemma \ref{lem-bdd-bel}, we obtain (\ref{new-T_3-inv}).


 Using the fact that $\Pi_h$ is a projection, and Lemma \ref{lem-meas}, (\ref{A-1}) and (\ref{new-pert-1}), we obtain,
\beq\label{new-pert-h-1} \|\Pi_hg_\gamma  - \Pi_h(g^{\e}_\gamma)\|_{L^2([0,1])}  \leq  \|g_\gamma - g^{\e}\circ \gamma\|_{L^2([0,1])} \leq \frac{\e}{C_{\gamma}},\eeq
and, using the fact that $\Pi_h$ is an orthogonal projection, and  (\ref{eq-proj-3}),
\begin{eqnarray}
\|g_\gamma - \Pi_hg_\gamma \|_{L^2([0,1])}
&=&  \sum_{\tau_h \in \mathbb{T}_h}\|g_\gamma - \Pi_hg_\gamma \|_{L^2(\tau_h)}
\leq  \sum_{\tau_h \in \mathbb{T}_h}\|g_\gamma \|_{L^2(\tau_h)}\nonumber \\
&\leq &  h^2(C_0)^4\sum_{\tau_h \in \mathbb{T}_h}\|g_\gamma\|_{H^4(\tau_h)}\leq 2(C_0)^4h^2\|g_\gamma\|_{H^4([0, 1])}\label{new-pert-h-2}  \end{eqnarray}
Taking $\tilde{C_0}=2(C_0)^4$, (\ref{new-pert-h-1}) and (\ref{new-pert-h-2}) imply     
\beq \label{new-pert-h} \|g_\gamma - \Pi_hg^{\e}_\gamma\|_{L^2([0,1])}
 \leq  h^2\tilde{C_0}\|g_\gamma\|_{H^4([0,1])}+ \frac{\e}{C_{\gamma}}. \eeq
Now, $\zeta \in \mathcal{W}$ implies $\zeta_{|_{I_h}} \in H^2(I_h)$.
Hence, taking $\phi_1$ and $\phi_2$ as ${g\circ\gamma}|_{\tau_h}$ and ${\Pi_hg^{\e}_{\gamma}}|_{\tau_h}$ respectively, in the first part of Lemma \ref{lem-interval-3}, (\ref{new-pert-h}) and (\ref{ineq-imbedding}), we have,
\beqarray
{\|T_3^{h,{\e}}(\zeta_{|_{\tilde{I}^h_{\e}}})-T_3\zeta\|}^2_{L^2(\cup_{\tau_h \in \mathbb{T}_h}[t^{h,\e}_0, t^{h,\e}_1])} &=& \sum_{\tau_h \in \mathbb{T}_h} {\|\zeta\circ g_\gamma-\zeta\circ\Pi_hg^{\e}_\gamma\|}^2_{L^2([t^{h,\e}_0, t^{h,\e}_1])}\\
& \leq & \Big(h^2\tilde{C_0}\|g_\gamma\|_{H^4([0,1])}+ \frac{\e}{C_{\gamma}}\Big)^2 \sum_{\tau_h \in \mathbb{T}_h}  {\|\zeta'\|}^2_{L^{\infty}{g(\gamma(\tau_h))}}\\
&\leq & (C_I)^2\Big(h^2\tilde{C_0}\|g_\gamma\|_{H^4([0,1])}+ \frac{\e}{C_{\gamma}}\Big)^2 \sum_{\tau_h \in \mathbb{T}_h}  {\|\zeta'\|}^2_{H^1 {(g(\gamma(\tau_h)))}}\\
&= & (C_I)^2\Big(h^2\tilde{C_0}\|g_\gamma\|_{H^4([0,1])}+ \frac{\e}{C_{\gamma}}\Big)^2 \|\zeta'\|^2_{H^1(I)}
\eeqarray
Hence,
\beq \label{new-eq-pert-0-1}\|T_3^{h,{\e}}(\zeta_{|_{\tilde{I}^h_{\e}}})-T_3\zeta\|_{L^2(\cup_{\tau_h \in \mathbb{T}_h}[t^{h,\e}_0, t^{h,\e}_1])} \leq C_I\Big(h^2\tilde{C_0}\|g\circ\gamma\|_{H^4([0,1])}+ \frac{\e}{C_{\gamma}}\Big) \|\zeta'\|_{H^1(I)}.\eeq
Since $g_\gamma'>0$, we have  $g(\gamma([t^{h,\e}_0, t^{h,\e}_1]))=[\tilde{g^{h,\e}}_0, \tilde{g^{h,\e}_1}] \subset I_h$ for some $\tilde{g^{h,\e}_0}< \tilde{g^{h,\e}_0}$. As $h$ and $\e$ satisfy (\ref{ineq-h}) and (\ref{A-h-2}), taking $\phi_1=(g\circ\gamma)|_{[t^{h,\e}_0, t^{h,\e}_1]}$ and $\phi_2=\Pi_h g^{\e}_{\gamma}|_{[t^{h,\e}_0, t^{h,\e}_1]}$ in Lemma \ref{lem-interval-1}, we have,
$$|g^{h,\e}_0- \tilde{g^{h,\e}_0}| < \tilde{C_0} h^{3/2} \|g_\gamma\|_{H^{4}(\tau_h)} +\frac{C'_0}{C_{\gamma}} \frac{\e}{h^{1/2}},$$
$$|g^{h,\e}_1- \tilde{g^{h,\e}_1}| < \tilde{C_0} h^{3/2} \|g_\gamma\|_{H^{4}(\tau_h)} +\frac{C'_0}{C_{\gamma}} \frac{\e}{h^{1/2}}.$$
Hence by (\ref{ineq-new-int}) and (\ref{ineq-new-int-1}),
\beq \label{new-ineq-interval-1} |g^h_0-\tilde{g^{h,\e}_0}|<2 \tilde{C_0} h^{3/2} \|g_\gamma\|_{H^{4}(\tau_h)} +\frac{2C'_0}{C_{\gamma}} \frac{\e}{h^{1/2}},\eeq
\beq \label{new-ineq-interval-2}|g^h_1-\tilde{g^{h,\e}_1}|<2\tilde{C_0} h^{3/2} \|g_\gamma\|_{H^{4}(\tau_h)} +\frac{2C'_0}{C_{\gamma}} \frac{\e}{h^{1/2}}.\eeq
Since (\ref{A-1}) and (\ref{A-2}) hold, by Lemma \ref{lem-meas},
\beqarray
{\|\zeta\circ g_\gamma \|}^2_{L^2([0,1]\setminus (\cup_{\tau_h \in \mathbb{T}_h}[t^{h,\e}_0, t^{h,\e}_1]) )} &=&
\int_{([0,1]\setminus (\cup_{\tau_h \in \mathbb{T}_h}[t^{h,\e}_0, t^{h,\e}_1]) )} {|\zeta(g(\gamma(s)))|}^2 ds \\
& \leq & \frac{2}{C_gC_{\gamma}} \int_{I\setminus \cup_{\tau_h \in \mathbb{T}_h} g(\gamma([t^{h,\e}_0, t^{h,\e}_1]))}
{|\zeta(z)|}^2 dz\\
& \leq & \frac{2}{C_gC_{\gamma}} \sum_{\tau_h \in \mathbb{T}_h} \int_{I_h\setminus g(\gamma([t^{h,\e}_0, t^{h,\e}_1]))} {|\zeta(z)|}^2 dz \\
& \leq & \frac{2}{C_gC_{\gamma}} \sum_{\tau_h \in \mathbb{T}_h} \left[ \int_{g^h_0}^{\tilde{g}^{h,\e}_0} {|\zeta(z)|}^2 dz + \int_{\tilde{g}^{h,\e}_1}^{g^h_1} |\zeta(z)|^2 dz \right].
\eeqarray
Hence,
\beq \label{new-eq-pert-1-1}{\|\zeta\circ g_\gamma \|}^2_{L^2([0,1]\setminus (\cup_{\tau_h \in \mathbb{T}_h}[t^{h,\e}_0, t^{h,\e}_1]) )}   \leq  \frac{2}{C_gC_{\gamma}} \sum_{\tau_h \in \mathbb{T}_h} \left[ \int_{g^h_0}^{\tilde{g}^{h,\e}_0} {|\zeta(z)|}^2 dz + \int_{\tilde{g}^{h,\e}_1}^{g^h_1} |\zeta(z)|^2 dz \right].
\eeq
Now, by (\ref{ineq-imbedding}), $\zeta \in \mathcal{W}$ implies $\zeta \in W^{1,\infty}(I)$. Hence,  as (\ref{new-ineq-interval-1}) and (\ref{new-ineq-interval-2}) hold, by Lemma \ref{lem-interval-2}-(i) and then by (\ref{ineq-imbedding}), we have
\beqarray
 \int_{g^h_0}^{\tilde{g}^{h,\e}_0} {|\zeta(z)|}^2 dz &\leq &  8(C_I)^2 \left(\tilde{C_0} h^{3/2} \|g_\gamma\|_{H^{4}(\tau_h)} +\frac{C'_0}{C_{\gamma}} \frac{\e}{h^{1/2}}\right)^2\|\zeta\|^2_{W^{1,\infty}([g^h_0, \tilde{g}^{h,\e}_0])}\\
 & \leq & 16(C_I)^4 \left(\tilde{C_0} h^{3/2} \|g_\gamma\|_{H^{4}(\tau_h)} +\frac{C'_0}{C_{\gamma}} \frac{\e}{h^{1/2}}\right)^2\|\zeta\|^2_{H^2([g^h_0, \tilde{g}^{h,\e}_0])},
 \eeqarray
and, similarly,
$$\int_{g^h_1}^{\tilde{g}^{h,\e}_1} {|\zeta(z)|}^2 dz \leq 16 (C_I)^4 \left(\tilde{C_0} h^{3/2} \|g_\gamma\|_{H^{4}(\tau_h)} +\frac{C'_0}{C_{\gamma}} \frac{\e}{h^{1/2}}\right)^2\|\zeta\|^2_{H^2([g^h_1, \tilde{g}^{h,\e}_1])}$$
Thus, from (\ref{new-eq-pert-1-1}) we have (\ref{new-eq-pert}).

Next, let $\zeta \in H^3(I)$. Since (\ref{new-ineq-interval-1}) and (\ref{new-ineq-interval-2}) hold, by Lemma \ref{lem-interval-2}-(ii) and then by (\ref{ineq-imbedding}), we obtain
\beqarray
 \int_{g^h_0}^{\tilde{g}^{h,\e}_0} {|\zeta(z)|}^2 dz &\leq & 32 (C_I)^4 \left(\tilde{C_0} h^{3/2} \|g_\gamma\|_{H^{4}(\tau_h)} +\frac{C'_0}{C_{\gamma}} \frac{\e}{h^{1/2}}\right)^3\|\zeta\|^2_{W^{2,\infty}([g^h_0, \tilde{g}^{h,\e}_0])} \\
 & \leq & 64 (C_I)^6 \left(\tilde{C_0} h^{3/2} \|g_\gamma\|_{H^{4}(\tau_h)} +\frac{C'_0}{C_{\gamma}} \frac{\e}{h^{1/2}}\right)^3\|\zeta\|^2_{H^3([g^h_0, \tilde{g}^{h,\e}_0])}, \eeqarray
 and, similarly,
 $$\int_{g^h_0}^{\tilde{g}^{h,\e}_0} {|\zeta(z)|}^2 dz \leq 64(C_I)^6  \left(\tilde{C_0} h^{3/2} \|g_\gamma\|_{H^{4}(\tau_h)} +\frac{C'_0}{C_{\gamma}} \frac{\e}{h^{1/2}}\right)^3\|\zeta\|^2_{H^3([g^h_1, \tilde{g}^{h,\e}_1])}. $$
Thus, from (\ref{new-eq-pert-1-1}) we have (\ref{new-eq-pert-1}).
\epf

Let $\tilde\zeta_{\e,\d,h} \in L^2(\tilde{I_{\e}})$ be the unique solution of the equation
\beq \label{new-eq-T3} {T_3^{h,\e}}^*T_3^{h,\e}(\zeta)={T_3^{h,\e}}^*f^{j^\d}, \eeq
that is  $\tilde\zeta_{\e,\d,h} := (T_3^{h,\e})^\dagger f^{j^\d}$.
Now, let $\zeta_{\e,\d,h} = \left\{\begin{array}{ll}
                            {\tilde\zeta}_{\e,\d,h}&\h{ on } \tilde{I_{\e}}\\
                            0& \h{ on }\, I\setminus \tilde{I_{\e}}.\end{array}\right.$
Then, $\zeta_{\e,\d,h}\in L^2(I)$.
Let $b_{\a,\e,\d,h}$ be the  solution of the equation
\beq \label{4.3.3}{(T^\a _2)}^* (T^\a _2) (w)= {(T^\a _2)}^*\zeta_{\e,\d,h}. \eeq
We show that $a_{\a,\e,\d,h}:=b'_{\a,\e,\d,h}$ is an approximate solution of (\ref{ill-1}). For this purpose, we shall make use of the following  proposition.

\begin{Proposition}
Let $a_0$ and $g$ be as defined in Lemma \ref{lem-err}. Let $h$ and $\e$ satisfy the relations in (\ref{ineq-h}) and (\ref{A-h-2}). Let $g^{\e} \in L^2(I)$ be such that (\ref{new-pert-1}) is satisfied. Then, $b_0=T_1(a_0)$ satisfies,
\beq \label{new-ineq-pf-2} \|b_0\|_{L^2(I\setminus \tilde{I_{\e}})}
  \leq \|b_0\|_{H^2(I)}(C_I)^2\Big(\tilde{C_0} h^{3/2} \|g_\gamma\|_{H^{4}([0,1])} +\frac{C'_0}{C_{\gamma}} \frac{\e}{h^{1/2}}\Big), \eeq
and, in addition, if $a_0\in  H^2(I)$, then,
\beq \label{new-ineq-pf-3} \|b_0\|_{L^2(I\setminus \tilde{I_{\e}})}
  \leq \|b_0\|_{H^3(I)}2(C_I)^3{\Big(\tilde{C_0} h^{3/2} \|g_\gamma\|_{H^{4}([0,1])} +\frac{C'_0}{C_{\gamma}} \frac{\e}{h^{1/2}}\Big)}^{3/2}, \eeq
 \end{Proposition}

\bpf
 Since, $h$ and $\e$ satisfy (\ref{ineq-h}), for any $\tau_h \in \mathbb{T}_h$, as (\ref{ineq-new-int}) holds, by Lemma \ref{lem-interval-2}-(i) and then by (\ref{ineq-imbedding}), we have
 \beqarray
 \|b_0\| _{L^2(I_h\setminus \tilde{I^h_{\e}})}
 &\leq & C_I \|b_0\|_{W^{1,\infty}(I_h)} \left(\tilde{C_0} h^{3/2} \|g_\gamma\|_{H^{4}(\tau_h)} +\frac{C'_0}{C_{\gamma}} \frac{\e}{h^{1/2}} \right) \\
 &\leq& (C_I)^2 \|b_0\|_{H^2(I_h)} \left(\tilde{C_0} h^{3/2} \|g_\gamma\|_{H^{4}(\tau_h)} +\frac{C'_0}{C_{\gamma}} \frac{\e}{h^{1/2}} \right),\eeqarray
and, if $a_0 \in H^2(I)$, $b_0 \in H^3(I)$ and so, by Lemma \ref{lem-interval-2}-(ii) and then by (\ref{ineq-imbedding}),
\beqarray
\|b_0\| _{L^2(I_h\setminus \tilde{I^h_{\e}})}
&\leq & 2(C_I)^2 \|b_0\|_{W^{2,\infty}(I_h)} \left(\tilde{C_0} h^{3/2} \|g_\gamma\|_{H^{4}(\tau_h)} +\frac{C'_0}{C_{\gamma}} \frac{\e}{h^{1/2}} \right)^{3/2}  \\
&\leq & 2(C_I)^3 \|b_0\|_{H^3(I_h)} \left(\tilde{C_0} h^{3/2} \|g_\gamma\|_{H^{4}(\tau_h)} +\frac{C'_0}{C_{\gamma}} \frac{\e}{h^{1/2}} \right)^{3/2}.\eeqarray
Since $\|b_0\| _{L^2(I\setminus \tilde{I_{\e}})}= \sum_{\tau_h\in \mathbb{T}_h} \|b_0\| _{L^2(I_h\setminus \tilde{I^h_{\e}})},$ the required inequalities follow.
\epf

\bt \label{new-Th-err}
Let $a_0$, $g$ and $j$ be as in Lemma \ref{lem-err}. Let $g^{\e} \in L^2(I)$, $j^{\d}\in W^{1-1/p,p}(\partial \O)$ with $p>3$. Also, let $g^{\e}$ and $j^{\d}$ satisfy (\ref{new-pert-1}) and (\ref{pert-2}), respectively, and $h$ and $\e$ satisfy the relations in  (\ref{A-h-2}) and (\ref{ineq-h}), and $a_{\a,\e,\d,h}=b'_{\a,\e,\d,h}$. Then the following results hold.
\vsq
{\rm(i)}\, With the original assumption that $a_0 \in H^1(I)$,
\begin{eqnarray}\label{new-1err-est-h1}  \|a_0-a_{\a,\e,\d,h} \|_{H^1(I)} &\leq & C_{\a} + \frac{2}{\a} \|b_0\|_{H^2(I)}(C_I)^2 C_{g,\e,h} \\
\nonumber && + \frac{2}{\a} \sqrt{C'_g C'_\gamma}   [D^1_{g,\e,h}\|b'_0\|_{H^1(I)} + D^2_{g,\e,h}\|b_0\|_{H^2(I)} + \tilde{C_{\gamma}}\delta],\\
\|a_0- a_{\a,\e,\d,h} \|_{L^2(I)} &\leq &\sqrt{\a}\|a'_0\|_{L^2(I)} +
\frac{2}{ \sqrt{\a}}   \|b_0\|_{H^2(I)}(C_I)^2 C_{g,\e,h} \\
\nonumber && + \frac{2}{ \sqrt{\a}}  \sqrt{C'_g C'_\gamma}   [D^1_{g,\e,h}\|b'_0\|_{H^1(I)} + D^2_{g,\e,h}\|b_0\|_{H^2(I)} + \tilde{C_{\gamma}}\delta].\end{eqnarray}

{\rm(ii)}\,
If $a_0\in H^2(I)$, then,
\begin{eqnarray} \label{new-1err-est-l2-2}
\|a_0- a_{\a,\e,\d,h} \|_{L^2(I)} &\leq & \sqrt{\a}\|a'_0\|_{L^2(I)} +
\frac{2}{ \sqrt{\a}} \|b_0\|_{H^3(I)} (C_I)^3 {(C_{g,\e,h})}^{3/2} \\
\nonumber && + \frac{2}{ \sqrt{\a}} \sqrt{C'_g C'_\gamma}   [D^1_{g,\e,h}\|b'_0\|_{H^1(I)} + D^3_{g,\e,h}\|b_0\|_{H^3(I)} + \tilde{C_{\gamma}}\delta]. \end{eqnarray}

{\rm(iii)}\,
If $a_0 \in H^3(I)$, then
\begin{eqnarray}  
\label{new-2err-est-h1}
\|a_0-a_{\a,\e,\d,h} \|_{H^1(I)} &\leq & (1+ C_L)\|a'_0\|_{H^2(I)}\a + \frac{2}{\a}  \|b_0\|_{H^3(I)} (C_I)^3 {(C_{g,\e,h})}^{3/2}\\ 
\nonumber && + \frac{2}{\a} \sqrt{C'_g C'_\gamma}   [D^1_{g,\e,h}\|b'_0\|_{H^1(I)} + D^3_{g,\e,h}\|b_0\|_{H^3(I)} + \tilde{C_{\gamma}}\delta],\\
\label{new-2err-est-l2} \|a_0- a_{\a,\e,\d,h} \|_{L^2(I)} &\leq & (1+ C_L)\|a'_0\|_{H^2(I)}\a +
\frac{2}{\sqrt{\a}} \|b_0\|_{H^3(I)} (C_I)^3 {(C_{g,\e,h})}^{3/2} \\
\nonumber && + \frac{2}{\sqrt{\a}} \sqrt{C'_g C'_\gamma}   [D^1_{g,\e,h}\|b'_0\|_{H^1(I)} + D^3_{g,\e,h}\|b_0\|_{H^3(I)} + \tilde{C_{\gamma}}\delta]. \end{eqnarray}
In the above  $C_{\a}>0$ is such that $C_{\a} \to 0$ as $\a\to 0$, 
$b_0=T_1(a_0)$, 
\beqarray
C_{g,\e,h} &=& (\tilde{C_0} h^{3/2} \|g_\gamma\|_{H^{4}([0,1])} +\frac{C'_0}{C_{\gamma}} \frac{\e}{h^{1/2}})\\
D^1_{g,\e,h} &=& C_I(\tilde{C_0}\|g\circ\gamma\|_{H^4([0,1])}h^2+ \frac{\e}{C_{\gamma}}), \\
D^2_{g,\e,h} &=&  4(C_I)^2\sqrt{\frac{2}{C_gC_{\gamma}}}{(\tilde{C_0} h^{3/2} \|g_\gamma\|_{H^{4}([0,1])} +\frac{C'_0}{C_{\gamma}} \frac{\e}{h^{1/2}})},\\
D^3_{g,\e,h} &=&  8(C_I)^3\sqrt{\frac{1}{C_gC_{\gamma}}} {(\tilde{C_0} h^{3/2} \|g_\gamma\|_{H^{4}([0,1])} +\frac{C'_0}{C_{\gamma}} \frac{\e}{h^{1/2}})}^{3/2},\eeqarray
and $C_0$, $C_L$, $\tilde{C_0}$, $C'_0$, $C_I$ $C_g$, $C'_g$, $C'_\gamma $, $C_\gamma$ are constants as defined in (\ref{eq-proj-1}), Proposition \ref{prop-projection}, Theorem \ref{th-new-pert-1}-(ii), (\ref{ineq-imbedding}), (\ref{A-2}) and (\ref{A-1}) respectively.
\et

\bpf
By definition of $\zeta_{\e,\d,h}$,
\beq \label{new-ineq-pf-1} \|\zeta_{\e,\d, h} - b_0\|_{L^2(I)} \leq \|\tilde{\zeta}_{\e,\d, h} - {b_0}|_{\tilde{I_{\e}}}\|_{L^2(\tilde{I_{\e}})} + \|b_0\|_{L^2(I\setminus \tilde{I_{\e}})}. \eeq
We use the notation  $(T_3^{h,\e})^\dagger := {({(T_3^{h,\e})}^*T_3^{h,\e})}^{-1}{(T_3^{h,\e})}^*$.
Then, by (\ref{new-T_3-inv}), and using the fact that $(T_3^{h,\e})^\dagger {(T_3^{h,\e})}^*$ is identity, we have
$$\|(T_3^{h,\e})^\dagger  T_3(T_2(b_0))-{(T_2(b_0))}|_{\tilde{I_{\e}}}\|_{L^2(\tilde{I_{\e}})} \leq \sqrt{2C'_g C'_\gamma} \|T_3(T_2(b_0))-T^{h,\e}_3({(T_2(b_0))}|_{\tilde{I_{\e}}})\|_{L^2([0,1])},$$
and, in addition, using (\ref{pert-3}),
\beq\label{new-th-er-1-2}\|(T_3^{h,\e})^\dagger (f^{j}-f^{j^\d})\|_{L^2(I)} \leq \sqrt{2C'_g C'_\gamma}\|f^{j^\d}-f^j\|_{L^2([0,1])}  \leq \sqrt{2C'_g C'_\gamma}\tilde{C_{\gamma}}\delta. \eeq
Also, since $a_0(g_1)=0$ we have $b_0=T_1(a_0) \in \mathcal{W}$, so that, by (\ref{new-eq-pert}) and (\ref{new-eq-pert-1}),
\beq \label{new-th-er-1-1} \|(T_3^{h,\e})^\dagger T_3(T_2(b_0))-{(T_2(b_0))}|_{\tilde{I_{\e}}}\|_{L^2(\tilde{I_{\e}})} \leq \sqrt{2C'_g C'_\gamma}   [D^1_{g,\e,h}\|b'_0\|_{H^1(I)} + D_{g,\e,h,b_0} ],\eeq
where,
$$D_{g,\e,h,b_0} := \left\{\begin{array}{ll}
                    D^2_{g,\e,h}\|b_0\|_{H^2(I)}&\h{  if }\, b_0\in \mathcal{W},\\
                    D^3_{g,\e,h}\|b_0\|_{H^3(I)} &\h{  if }\, b_0\in H^3(I) \cap \mathcal{W}.\end{array}\right.$$
Now, by the definition of $\tilde{\zeta}_{\e,\d, h}$ and the fact that $T_3(T_2(b_0))=f^j$, we have
\beqarray
\|\tilde{\zeta}_{\e,\d, h} - {b_0}|_{\tilde{I_{\e}}}\|_{L^2(\tilde{I_{\e}})}
 &\leq & \|(T_3^{h,\e})^\dagger  f^{j^\d}-{(T_2(b_0))}|_{\tilde{I_{\e}}}\|_{L^2(\tilde{I_{\e}})} \\
 &\leq & \|(T_3^{h,\e})^\dagger  T_3(T_2(b_0))-{(T_2(b_0))}|_{\tilde{I_{\e}}}\|_{L^2(\tilde{I_{\e}})} \\
 && +
 \|(T_3^{h,\e})^\dagger (f^{j}-f^{j^\d})\|_{L^2(I)}
 \eeqarray
Hence, from (\ref{new-th-er-1-1}) and (\ref{new-th-er-1-2}) we have
 \beq \label{new-th-er-1-3} \|\tilde{\zeta}_{\e,\d, h} - {b_0}|_{\tilde{I_{\e}}}\|_{L^2(\tilde{I_{\e}})} \leq
\sqrt{2C'_g C'_\gamma}   [D^1_{g,\e,h}\|b_0'\|_{H^1(I)} + D_{g,\e,h,b_0} + \tilde{C_{\gamma}}\delta].\eeq
Thus, from (\ref{new-ineq-pf-1}), (\ref{new-ineq-pf-2}) and (\ref{new-th-er-1-3}) we have
\begin{eqnarray}
\label{new-th-er-1-4} \|\zeta_{\e,\d, h} - b_0\|_{L^2(I)}
&\leq &  \|b_0\|_{H^2(I)}(C_I)^2\Big(\tilde{C_0}\|g\circ\gamma\|_{H^4([0,1])}h^2 + \frac{\e}{C_{\gamma}}\Big) \\
\nonumber && + \sqrt{2C'_g C'_\gamma}   [D^1_{g,\e,h}\|b_0'\|_{H^1(I)} + D^2_{g,\e,h}\|b_0\|_{H^2(I)} + \tilde{C_{\gamma}}\delta].\end{eqnarray}
If $a_0 \in H^2(I)$ then $b_0 \in H^3(I)$, and thus from (\ref{new-ineq-pf-1}), (\ref{new-ineq-pf-3}) and (\ref{new-th-er-1-3}) we have,
\begin{eqnarray}
\label{new-th-er-1-5}\|\zeta_{\e,\d, h} - b_0\|_{L^2(I)}
&\leq & \|b_0\|_{H^3(I)}(C_I)^3{\left(\tilde{C_0}\|g\circ\gamma\|_{H^4([0,1])}h^2+ \frac{\e}{C_{\gamma}}\right)}^{3/2}\\
\nonumber && + \sqrt{2C'_g C'_\gamma}   [D^1_{g,\e,h}\|b_0'\|_{H^1(I)} + D^3_{g,\e,h}\|b_0\|_{H^3(I)} + \tilde{C_{\gamma}}\delta]. \end{eqnarray}
Our aim is to find an estimate for the error term $(a_0-a_{\a,\e,\d,h})$ in $L^2(I)$ and $H^1(I)$ norms. Now $b_{\a,\e,\d,h}$ is the unique solution of equation (\ref{4.3.3}). Thus, according to Lemma \ref{lem-err} we need an estimate of $\|\zeta_{\e,\d, h} - b_0\|_{L^2(I)}$ in order to find our required estimates.
Inequalities (\ref{new-th-er-1-4}) and (\ref{new-th-er-1-5}) give us estimates of $\|\zeta_{\e,\d, h} - b_0\|_{L^2(I)}$ under different conditions on $b_0$. Hence, taking $\zeta_{\e,\d,h}$ in place of $\zeta$ in Lemma \ref{lem-err} we have the proof.
\epf

\brem \label{rem-order-1}
Suppose
$$2{\e}^{1/4}< \min\left\{ \left({\frac{C_{\gamma}C_g}{\tilde{C}_1\|g\circ\gamma\|_{H^4([0,1])}+\frac{C'_1}{C'_{\gamma}}}}\right)^{1/2},\q  \frac{1}{\tilde{C_0}\|g\circ\gamma\|_{H^4([0,1])}+\frac{C'_0}{C'_{\gamma}}} \right\}.$$
Then, for  $\e=\d$ and $h={\d}^{1/2}$,  (\ref{A-h-2}) and (\ref{ineq-h}) are satisfied. Hence, by Theorem \ref{new-Th-err}, we have the following:

\ben
\i Choosing $\a=\sqrt{\d}$, we have
$$\|a_0- a_{\a,\e,\d,h}\|_{H^1(I)}=o(1).$$
\i If $a_0 \in H^3(I)$ and $\a={\d}^{2/3}$, then
$$\|a_0- a_{\a,\e,\d,h}\|_{H^1(I)}=O(\sqrt{\d}),\q \|a_0- a_{\a,\e,\d,h}\|_{L^2(I)}= O({\d}^{2/3}).$$
\i Choosing $\a=\d$, we have
$$\|a_0- a_{\a,\e,\d,h}\|_{L^2(I)}=O(\d^{1/4}).$$
\i
If $a_0\in H^2(I)$, then  $$\|a_0- a_{\a,\e,\d,h}\|_{L^2(I)}= O({\d}^{1/2}).$$\een
Results in (1) and (2) above  are analogous to the corresponding results for $a_0- a_{\a,\e,\d}$  in Remark \ref{rem-err-est-1}.
The  estimate in (4) is same as the corresponding estimate in Remark \ref{rem-err-est-1}, except for the fact that here we need an additional condition that $a_0 \in H^2(I)$.
\erem

\section{With exact solution having non-zero value at $g_1$} \label{sec-5}
In the previous two sections we have considered the exact solution with assumption that $a_0(g_1)=0$. Here we consider the case when $a_0(g_1)\neq 0$ but is assumed to be known.
Let $a_0(g_1)=c$.
Since $a_0$ is the solution to Problem (P), by (\ref{ill-1}) we have
$f^j=T(a_0)$ which implies
\beq \label{5.1} f^j = T(a_0-c+c)=T(a_0-c)+cT(1)\eeq
Now by definition of $T$ we have
 \beq \label{5.2} T(1)(s)=\int_{g_0}^{g\circ\gamma(s)}dt=g\circ\gamma(s)-g_0, \q s\in [0,1].\eeq
Thus, combining (\ref{5.1}) and (\ref{5.2}) we have
\beq \label{5.3} T(a_0-c)= f^j -c(g_\gamma - g_0)\eeq
Hence $a_0-c$ is the solution of the following operator equation,
\beq \label{new-op-eq} T(a)=f^j -c(g_\gamma - g_0),\eeq
where clearly $f^j -c(g_\gamma - g_0) \in L^2([0,1])$.
Also, $(a_0-c)(g_1)=0$. Now, let us define
$$b_{0,c}(x)=\int_{g_0}^x (a_0(t)-c)dt, \q x\in I.$$
Then $b_{0,c} \in \mathcal{W}$.
Thus, the analysis of the previous two sections can be applied here to obtain a stable approximate solution of equation (\ref{new-op-eq}).
Let $a_{c,\a}:=b'_{c,\a}$, where $b_{c,\a}$ is the solution to the following equation.
\beq \label{op-eq-1} {(T^\a _2)}^* (T^\a _2) (w)= {(T^\a _2)}^*\zeta_c, \eeq
where $\zeta_c$ is the solution of the equation
\beq \label{op-eq-2} {(T_3)}^*(T_3)\zeta= {(T_3)}^*(f^j -c(g_\gamma-g_0)).\eeq
Now, let $g^{\e}$ and $j^{\d}$ be the perturbed data as defined in Theorem \ref{new-Th-err}. Also, let $g$ be such that $g\circ \g \in H^4([0,1])$.
Let $\tilde{\zeta}_{c,\e,\d,h}$ be the solution of the equation
\beq \label{5.3.2}{T_3^{h,\e}}^*T_3^{h,\e}(\zeta)={T_3^{h,\e}}^*(f^{j^{\d}} -c(\Pi_h(g^{\e}_{\gamma})-g_0)),\eeq
where $\Pi_h(g^{\e}_{\gamma})$ is as defined in Section \ref{sec-4}.
Now, $\zeta_{c,\e,\d,h}$ defined as, $\tilde{\zeta}_{c,\e,\d,h}$ on $\tilde{I_{\e}}$ and $0$ on $I\setminus \tilde{I_{\e}}$, is in $L^2(I)$.
Let $b_{c,\e,\d,h}$ be the solution of the equation
\beq \label{5.3.3} {(T^\a _2)}^* (T^\a _2) (w)= {(T^\a _2)}^*\zeta_{c,\e,\d,h}\eeq
Then we have the following theorem.

\bt \label{6-th-err}
Let $a_0$, $c$ and $b_{0,c}$ be as defined in the beginning of the section. Let $g$ and $j$ be as defined in Lemma \ref{lem-err}, and $g\circ\gamma \in H^4([0,1])$. Let $h$ and $\e$ satisfy (\ref{ineq-h}) and (\ref{A-h-2}), respectively. Also, let $g^{\e} \in L^2(\Gamma)$, $j^{\d}\in W^{1-1/p,p}(\partial \O)$ with $p>3$, and $g^{\e}$ and $j^{\d}$ satisfy (\ref{new-pert-1}) and (\ref{pert-2}) respectively. Let $a_{c,\a,\e,\d,h}:=b'_{c,\a,\e,\d,h}$, and let 
\beqarray
C_{g,\e,h} &:=&  \tilde{C_0} h^{3/2} \|g_\gamma\|_{H^{4}(\tau_h)} +\frac{C'_0}{C_{\gamma}} \frac{\e}{h^{1/2}},\\
D^1_{g,\e,h} &=&  C_I(\tilde{C_0}\|g\circ\gamma\|_{H^4([0,1])}h^2+ \frac{\e}{C_{\gamma}}),\\
D_{g,\e,h,b_{0,c}} &=&  \left\{\begin{array}{ll}
                    4(C_I)^2\sqrt{\frac{2}{C_gC_{\gamma}}}{(\tilde{C_0} h^{3/2} \|g_\gamma\|_{H^{4}(\tau_h)} +\frac{C'_0}{C_{\gamma}} \frac{\e}{h^{1/2}})} &\h{ if }\, b_{0,c}\in \mathcal{W},\\
                     8(C_I)^3\sqrt{\frac{1}{C_gC_{\gamma}}} {\left(\tilde{C_0} h^{3/2} \|g_\gamma\|_{H^{4}(\tau_h)} +\frac{C'_0}{C_{\gamma}} \frac{\e}{h^{1/2}}\right)}^{3/2} &\h{  if }\, b_{0,c}\in H^3(I) \cap \mathcal{W}.\end{array}\right.\\
\mathcal{H}(c,\e,\d,h) &:=&   \sqrt{C'_gC'_{\gamma}} \left(D^1_{g,\e,h}\|b'_{0,c}\|_{H^1(I)} + D_{g,\e,h,b_{0,c}} + \tilde{C_{\gamma}}\d + c D^1_{g,\e,h}\right).
\eeqarray
Then
\begin{eqnarray} \label{5-1err-est-h1}
\|a_0- (c+a_{c,\a,\e,\d,h})\|_{H^1(I)} &\leq&  C_{\a} + \frac{2}{\a}\left( \|b_{0,c}\|_{H^2(I)}C_I^2 C_{g,\e,h}+\mathcal{H}(c,\e,\d,h)\right),\\
\label{5-1err-est-l2} \|a_0-(c+a_{c,\a,\e,\d,h})\|_{L^2(I)} &\leq& \sqrt{\a}\|a'_0\|_{L^2(I)} +
\frac{2}{\sqrt{\a}}\left( \|b_{0,c}\|_{H^2(I)}C_I^2 C_{g,\e,h}+\mathcal{H}(c,\e,\d,h) \right), 
\end{eqnarray}
where $C_{\a}>0$ is such that $C_{\a} \to 0$ as $\a\to 0$.
Further,  we have the following.
\ben
\i[\rm(i)] 
If $a'_0 \in L^{\infty}(I)$, then,
\beq \label{new-1err-est-l2-2} \|a_0- (c+ a_{\a,\e,\d,h}) \|_{L^2(I)} \leq \sqrt{\a}\|a'_0\|_{L^2(I)} +
\frac{2}{ \sqrt{\a}}\left(\|b_{0,c}\|_{H^3(I)} (C_I)^3 {(C_{g,\e,h})}^{3/2} + \mathcal{H}(c,\e,\d,h)\right). \eeq
\i[\rm(ii)] If $a_0 \in H^3(I)$, then
\beq \label{5-2err-est-h1} \|a_0- (c+a_{c,\a,\e,\d,h})\|_{H^1(I)} \leq (1+ C_L)\|a'_0\|_{H^2(I)}\a + \frac{2}{\a}
\left( \|b_{0,c}\|_{H^3(I)}(C_I)^3 {(C_{g,\e,h})}^{3/2} + \mathcal{H}(c,\e,\d,h)\right), \eeq
\beq \label{5-2err-est-l2} \|(a_0-c)-a_{c,\a,\e,\d,h}\|_{L^2(I)} \leq (1+ C_L)\|a'_0\|_{H^2(I)}\a +
\frac{2}{\sqrt{\a}} \left( \|b_{0,c}\|_{H^3(I)}(C_I)^3 {(C_{g,\e,h})}^{3/2} + \mathcal{H}(c,\e,\d,h)\right),\eeq
 \een
where $C_L$, $\tilde{C_0}$, $C'_0$, $C_I$ $C_g$, $C'_g$, $C'_\gamma $, and $C_\gamma$ are constants as in Proposition \ref{prop-projection}, Theorem \ref{th-new-pert-1},  (\ref{ineq-imbedding}), (\ref{A-2}) and (\ref{A-1}) respectively.
\et

\bpf
By definition of $\zeta_{c,\e,\d,h}$,
\beq \label{6-ineq-pf-1} \|\zeta_{c,\e,\d, h} - b_{0,c}\|_{L^2(I)} \leq \|\tilde{\zeta}_{c,\e,\d, h} - {b_{0,c}}|_{\tilde{I_{\e}}}\|_{L^2(\tilde{I_{\e}})} + \|b_{0,c}\|_{L^2(I\setminus \tilde{I_{\e}})}. \eeq
Here also we use the notation $(T_3^{h,\e})^\dagger := {({(T_3^{h,\e})}^*T_3^{h,\e})}^{-1}{(T_3^{h,\e})}^*$.
By (\ref{new-T_3-inv}),
$$\|(T_3^{h,\e})^\dagger T_3(T_2(b_{0,c}))-{(T_2(b_{0,c}))}|_{\tilde{I_{\e}}}\|_{L^2(\tilde{I_{\e}})} \leq \sqrt{2C'_g C'_\gamma} \|T_3(T_2(b_{0,c}))-T^{h,\e}_3({(T_2(b_{0,c}))}|_{\tilde{I_{\e}}})\|_{L^2([0,1])},$$
\beqarray
\|(T_3^{h,\e})^\dagger (f^{j}-c(g_\gamma-g_0)- f^{j^\d}-c(\Pi_hg^{\e}_\gamma - g_0))\|_{L^2(I)} 
 &\leq &  \sqrt{2C'_g C'_\gamma} \|[f^{j}-c(g_\gamma-g_0)\\
 &&- f^{j^\d}- c(\Pi_hg^{\e}_\gamma  -g_0)]\|_{L^2([0,1])}\\
&\leq &   \sqrt{2C'_g C'_\gamma} \|[f^{j}-f^{j^\d}\\
&&- c(g\circ\gamma-\Pi_hg^{\e}_\gamma)] \|_{L^2(\tilde{I_{\e}})}.
\eeqarray
Hence,
 by (\ref{pert-3}) and (\ref{new-pert-h}),
\beq \label{6-th-er-1-2} \|(T_3^{h,\e})^\dagger (f^{j}-c(g_\gamma-g_0)-f^{j^\d}-c(\Pi_hg^{\e}_\gamma  -g_0))\|_{L^2(I)} \leq \sqrt{2C'_g C'_\gamma}(\tilde{C_{\gamma}}\d + c D^1_{g,\e,h}), \eeq
and,
by definition $b_{0,c} \in \mathcal{W}$ and so, by (\ref{new-eq-pert}) and (\ref{new-eq-pert-1})
\beq \label{6-th-er-1-1} \|(T_3^{h,\e})^\dagger T_3(T_2(b_{0,c}))-{(T_2(b_{0,c}))}|_{\tilde{I_{\e}}}\|_{L^2(\tilde{I_{\e}})} \leq \sqrt{2C'_g C'_\gamma}   [D^1_{g,\e,h}\|b'_{0,c}\|_{H^1(I)} + D_{g,\e,h,b_{0,c}} ].\eeq
Now by definition of $\tilde{\zeta}_{c,\e,\d, h}$ and the fact that $T_3(T_2(b_{0,c}))=f^j-c(g_\gamma-g_0)$, we have
\beqarray
\|\tilde{\zeta}_{c,\e,\d, h} - {b_{0,c}}|_{\tilde{I_{\e}}}\|_{L^2(\tilde{I_{\e}})}
 &\leq & \|(T_3^{h,\e})^\dagger (f^{j^\d}-c(\Pi_hg^{\e}_\gamma -{(T_2(b_{0,c}))}|_{\tilde{I_{\e}}}\|_{L^2(\tilde{I_{\e}})} \\
 &\leq & \|(T_3^{h,\e})^\dagger (T_3(T_2(b_{0,c}))-{(T_2(b_{0,c}))}|_{\tilde{I_{\e}}})\|_{L^2(\tilde{I_{\e}})} \\
 && +
 \|(T_3^{h,\e})^\dagger (f^{j}-c(g_\gamma-g_0)-f^{j^\d}-c(\Pi_hg^{\e}_\gamma  -g_0))\|_{L^2(I)}.
 \eeqarray
Hence, from (\ref{6-th-er-1-1}) and (\ref{6-th-er-1-2}) we have
 \beq \label{6-th-er-1-3} \|\tilde{\zeta}_{\e,\d, h} - {b_{0,c}}|_{\tilde{I_{\e}}}\|_{L^2(\tilde{I_{\e}})} \leq
\sqrt{2C'_g C'_\gamma}   [D^1_{g,\e,h}\|b'_{0,c}\|_{H^1(I)} + D_{g,\e,h,b_{0,c}} + \tilde{C_{\gamma}}\d + c D^1_{g,\e,h}].\eeq
Thus, from (\ref{6-ineq-pf-1}), (\ref{new-ineq-pf-2}) and (\ref{6-th-er-1-3}) we have
\beq \label{6-th-er-1-4}\|\zeta_{\e,\d, h} - b_{0,c}\|_{L^2(I)} \leq \|b_{0,c}\|_{H^2(I)}(C_I)^2\left(\tilde{C_0}\|g\circ\gamma\|_{H^4([0,1])}h^2+ \frac{\e}{C_{\gamma}}\right) + \mathcal{H}(c,\e,\d,h).\eeq
If $a_{0,c} \in H^2(I)$, from (\ref{6-ineq-pf-1}), (\ref{new-ineq-pf-3}) and (\ref{6-th-er-1-3}) we have,
\beq \label{new-th-er-1-5}\|\zeta_{\e,\d, h} - b_{0,c}\|_{L^2(I)} \leq \|b_{0,c}\|_{H^3(I)}(C_I)^3{\left(\tilde{C_0}\|g\circ\gamma\|_{H^4([0,1])}h^2+ \frac{\e}{C_{\gamma}}\right)}^{3/2} + \mathcal{H}(c,\e,\d,h).\eeq
By definition, $b_{c,\a,\e,\d,h}$ is the unique solution of equation (\ref{5.3.3}). 
Thus, putting $\zeta_{c,\e,\d,h}$ in place of $\zeta$ in Lemma \ref{lem-err}, we have the proof using (\ref{new-th-er-1-4}) and (\ref{new-th-er-1-5}).
\epf

From Theorem \ref{6-th-err}, we see that $c+ a_{c,\a,\e,\d,h}$ is a stable approximate solution of Problem (P), with error estimates obtained from Theorem \ref{6-th-err}.

\brem
Let us relax the assumption on the exact solution $a_0$ even more. Let us assume that $a_0(g_1)$ is not equal to the known number $c$ but is known to be ``close" to it, i.e,
\beq \label{5-approx} |a_0(g_1)-c| < \eta,\eeq
for some $\eta>0$.
Let $c_0:=a_0(g_1)$.
Define $b_{0,c_0}(x)= \int_{g_0}^x (a_0(t)-c_0)dt$ for $x\in I$. Then $b_{0,c_0}\in \mathcal{W}$. Also, let $g$, $j$, $g^{\e}$, $j^{\d}$, $h$, $\zeta_{c,\e,\d,h}$, $b_{c,\a,\e,\d,h}$ and $a_{c,\a,\e,\d,h}$ be as defined in
Theorem \ref{6-th-err}.
Since (\ref{5-approx}) holds,
\beqarray
\|(f^{j}-f^{j^\d}-(c-c_0)(g_\gamma+g_0)
   - c(g\circ\gamma-\Pi_hg^{\e}_\gamma )\|_{L^2(\tilde{I_{\e}})}\|_{L^2([0,1])}& \leq &  \|f^{j^\d}-f^j\|_{L^2([0,1])}\\
   & + & c \|g\circ\gamma-\Pi_hg^{\e}_\gamma \|_{L^2([0,1])} \\
   &+& (\|g\circ\gamma\|_{L^2([0,1])}+|g_0|)\eta
   \eeqarray
 and, by (\ref{pert-3}) and (\ref{new-pert-h}), we have
 \beq \label{7-th-1} \|(f^{j}-f^{j^\d}-(c-c_0)(g_\gamma+g_0)
   - c(g\circ\gamma-\Pi_hg^{\e}_\gamma )\|_{L^2(\tilde{I_{\e}})}\|_{L^2([0,1])} \leq \tilde{C_
 {\gamma}}\d +\, a_0(g_1)D^1_{g,\e,h}+(\|g\circ\gamma\|_{L^2([0,1])}+|g_0|)\eta, \eeq
 with $D^1_{g,\e,h}$ as in Theorem \ref{6-th-err}.
 Now, as $T_3(T_2(b_{0,c_0}))=f^j-c_0(g_\gamma-g_0)$,
\beqarray
\|\tilde{\zeta}_{c,\e,\d,h} - b_{0,c_0}|_{L^2(\tilde{I_{\e}})}\|_{L^2(\tilde{I_{\e}})} &=& \|\tilde{\zeta}_{c,\e,\d,h} - T_2(b_{0,c_0})|_{L^2(\tilde{I_{\e}})}\|_{L^2(\tilde{I_{\e}})}\\
& \leq & \|(T_3^{h,\e})^\dagger (f^{j^\d}-c(\Pi_hg^{\e}_\gamma  -g_0))-T_2(b_{0,c_0})|_{L^2(\tilde{I_{\e}})}\|_{L^2(\tilde{I_{\e}})} \\
&\leq & \|(T_3^{h,\e})^\dagger T_3(T_2(b_{0,c}))-T_2(b_{0,c})|_{L^2(\tilde{I_{\e}})}\|_{L^2(\tilde{I_{\e}})}\\
 && +    \|(T_3^{h,\e})^\dagger (f^{j}-c_0(g_\gamma-g_0)-f^{j^\d}+
 c(\Pi_hg^{\e}_\gamma  -g_0))\|_{L^2(\tilde{I_{\e}})}\\
   &\leq & \|(T_3^{h,\e})^\dagger (T_3T_2(b_{0,c})-T_3^{h,\e}(T_2(b_{0,c})|_{L^2(\tilde{I_{\e}})})\|_{L^2(\tilde{I_{\e}})}\\
 && + \|(T_3^{h,\e})^\dagger (f^{j}-f^{j^\d}-(c-c_0)(g_\gamma+g_0)
   - c(g\circ\gamma-\Pi_hg^{\e}_\gamma ))\|_{L^2(\tilde{I_{\e}})},
\eeqarray
and, by (\ref{new-T_3-inv}), (\ref{new-eq-pert}) and (\ref{new-eq-pert-1})
\beqarray \|(T_3^{h,\e})^\dagger T_3(T_2(b_{0,c}))-{(T_2(b_{0,c}))}|_{\tilde{I_{\e}}}\|_{L^2(\tilde{I_{\e}})}
& \leq &  \sqrt{2C'_g C'_\gamma}   [D^1_{g,\e,h}\|b'_{0,c}\|_{H^1(I)} +   D_{g,\e,h,b_{0,c}} \\
&&+ \| f^{j}-f^{j^\d} \|_{L^2([0,1])}\\
&&+ \|(c-c_0)(g_\gamma+g_0)+c(g\circ\gamma-\Pi_hg^{\e}_\gamma )\|_{L^2([0,1])} ]\eeqarray
with $D_{g,\e,h,b_{0,c}}$ as in Theorem \ref{6-th-err}.
Hence, by (\ref{7-th-1})
\beqarray
\|\tilde{\zeta}_{c,\e,\d,h} - b_{0,c_0}|_{L^2(\tilde{I_{\e}})}\|_{L^2(\tilde{I_{\e}})} &\leq & \sqrt{2C'_g C'_\gamma}   [D^1_{c,g,\e,h}\|b'_{0,c}\|_{H^1(I)} + D_{g,\e,h,b_{0,c}} + \tilde{C_
 {\gamma}}\d \\
 &&+\, cD^1_{g,\e,h}+(\|g\circ\gamma\|_{L^2([0,1])}+|g_0|)\eta ].
\eeqarray
Thus, using similar arguments as  in the proof of Theorem \ref{6-th-err}, we obtain estimates for  
$$\|(a_0-c_0)-a_{c,\a,\e,\d,h}\|_{H^1(I)}\q\h{and}\q \|(a_0-c_0)-a_{c,\a,\e,\d,h}\|_{L^2(I)}.$$
Using the fact that 
$$(a_0-(a_{c,\a,\e,\d,h}+c))=((a_0-c_0)-a_{c,\a,\e,\d,h})+(c_0-c),$$ 
we obtain $(a_{c,\a,\e,\d,h}+c)$ as a stable approximate solution to Problem (P), and obtain the corresponding error estimates.
\erem

\section{Illustration of the procedure} \label{sec-6}
In order to find a stable approximate solution of Problem (P) using the new regularization method we have to undertake the following.

Let $j^{\d} \in W^{1-1/p,p}(\partial \O)$ with $p>3$, $g^{\e}\in L^2(\partial \O)$ be the perturbed data satisfying (\ref{pert-2}) and (\ref{new-pert-1}) respectively, and let
$f^{j^\d}=v^{j^{\d}}\circ \gamma$.
Also let us assume $g\circ \g \in H^4([0,1])$.
Then,  by the following steps we obtain the regularized solution $a_{\a,\e,\d}$.
\vsq

\begin{itemize}
\i[\underline{Step\,(i)}:] (a)\, Suppose $g^{\e} \in W^{1,\infty}(\Gamma)$ and it satisfies (\ref{pert-1}). Let $\tilde{\zeta} _{\e,\d}$  be the unique element in $L^2([0, 1])$ such that
\beq \label{ineq-7-1}{(T^{\e}_3)}^*(T^{\e}_3) \tilde{\zeta} _{\e,\d} = {(T^{\e}_3)}^*f^{j^{\delta}} \eeq  with
$T^{\e}_3$ defined as in (\ref{op-defn-1}).
Define $\zeta_{\e,\d}$ to be equal to $\tilde{\zeta}_{\e,\d}$ on $\tilde{I}_{\e}$, and equal to $0$ on $I\setminus \tilde{I}_{\e}$.
\vsq
\noi
(b)\, Suppose $g^{\e}\in L^2(\Gamma)\setminus W^{1,\infty}(\Gamma)$. Then under the assumption $g\circ \gamma\in H^4([0,1])$, there exists a unique element $\tilde{\zeta}_{\e,\d,h}\in L^2([0, 1])$ such that
\beq \label{ineq-7-2}{(T^{h,\e}_3)}^*(T^{h,\e}_3) \tilde{\zeta}_{\e,\d,h} = {(T^{h,\e}_3)}^*f^{j^{\delta}}\eeq with
$T^{h,\e}_3$ defined as in (\ref{op-defn-2}). Define $\zeta_{\e,\d,h}$ to be equal to $\tilde{\zeta}_{\e,\d,h}$ on $\tilde{I}_{\e}$, and equal to $0$ on $I\setminus \tilde{I}_{\e}$.

 We denote the solution obtained in this step by $\zeta^{\e,\d}$.
\vsq 
\i[\underline{Step\,(ii)}:] Let $\zeta^{\e,\d}$ be as in Step (i).
Let  $b_{\a,\e,\d}$ be the unique element in $H^2(I)$ such that 
\beq \label{ineq-7-3}(T_2^\a)^*T_2^\a (b_{\a,\e,\d}) = (T_2^\a)^*\zeta^{\e,\d}\eeq
with $T_2^\a$ defined as in (\ref{def-reg-T2}).

\vsq
\i[\underline{Step\,(iii)}:] Define $a_{\a,\e,\d}:= b'_{\a,\e,\d}$, the derivative of $b_{\a,\e,\d}$.
\end{itemize}
\vsh

We now explain how to solve (\ref{ineq-7-1}) and (\ref{ineq-7-2}) and obtain $\zeta^{\e,\d}$.
Let us observe that, for $g^{\e}\in W^{1,\infty}(\Gamma)$ and  for  $f \in L^2([0,1])$,
$$(T_3^{\epsilon})^*(f)(z)=   \left\{
\begin{array}{ll}
       \frac{f(( g^{\e} \circ \g)^{-1}(z))}{(g^{\e}\circ \gamma)'(\gamma ^{-1}( (g^{\e})^{-1}(z)))}
         & z\in \tilde{I}_\e  \\
       0 & z \in I\setminus \tilde{I}_\e.
\end{array}
\right. $$
Hence, it can be seen that,
$$\zeta_{\e,\d}(z)=   \left\{
\begin{array}{ll}
       f^{j^\d}(( g^{\e}\circ \g)^{-1}(z))
         & z\in \tilde{I}_\e  \\
       0 & z \in I\setminus \tilde{I}_\e.
\end{array}
\right. $$

For $g^{\e}\in L^2(\G)\setminus W^{1,\infty}(\G)$, for any $f \in L^2([0,1])$
$$(T^{h,\e}_3)^*(f)(z):= (S^{h,\e,\d})^*(f)(z), \q z\in \tilde{I}^h_\e,$$
where 
$$(S^{h,\epsilon,\d})^*(f)(z)=   \left\{
\begin{array}{ll}
       \frac{f((\Pi_hg^{\e} \circ \g)^{-1}(z))}{(\Pi_hg^{\e}\circ \gamma)'(\gamma ^{-1}( (\Pi_hg^{\e})^{-1}(z)))}
         & z\in \tilde{I}^h_\e  \\
       0 & z \in I_h\setminus \tilde{I}^h_\e,
\end{array}
\right. $$
Hence, it can be seen that,  
$$\zeta_{h,\e,\d}(z)= \chi^{h,\e,\d}(z), \q z\in \tilde{I}^h_\e, $$
where,
$$\chi_{h,\e,\d}(z)=   \left\{
\begin{array}{ll}
       f^{j^\d}((\Pi_hg^{\e}\circ \g)^{-1}(z))
         & z\in \tilde{I}^h_\e  \\
       0 & z \in I_h\setminus \tilde{I}^h_\e.
\end{array}
\right. $$
Thus we have $\zeta^{\e,\d}$.
Next let us consider Step (ii).
Let us consider the case when $\zeta^{\e,\d} \in C(I)$. If $\zeta^{\e,\d} \in R(T^{\a}_2)$ then the solution of 
\beq \label{ineq-7-4} T^{\a}_2 (b)= \zeta^{\e,\d}\eeq
is the solution of (\ref{ineq-7-3}). Now let us note that, finding a solution of (\ref{ineq-7-4}) is same as solving the ODE
\beq \label{ODE-1}-\a b''+ b=\zeta^{\e,\d} \eeq
with boundary condition
\beq \label{ODE-2} b(g_0)=0 \eeq
and
\beq \label{ODE-3} b'(g_1)=0. \eeq
Hence, if $j^{\d}$ and $g^{\e}$ are such that $\zeta^{\e,\d} \in R(T^{\a}_2)\cap C(I)$ then the solution of the ODE (\ref{ODE-1})-(\ref{ODE-3}) gives us our desired $b_{\a,\e,\d}$. Also, by Step (iii) $a_{\a,\e,\d}=b'_{\a,\e,\d}$ is our desired regularized solution. Now let us note that, if $\zeta^{\e,\d} \in L^2(I)\setminus C(I)$ then there exists $\zeta^{\e,\d}_n \in C(I)$ such that $$\|\zeta^{\e,\d}-\zeta^{\e,\d}_n\|_{L^2(I)}=O(\frac{1}{n})$$ for $n\in \N$. Since by (\ref{bdd-below-1}) we have
$$\|((T_2^\a)^*T_2^\a)^{-1}(T_2^\a)^*(\zeta^{\e,\d}-\zeta^{\e,\d}_n) \|_{H^2(I)} \leq \frac{1}{\a}\|\zeta^{\e,\d}-\zeta^{\e,\d}_n\|_{L^2(I)},$$
if $\zeta^{\e,\d}_n \in R(T^{\a}_2)$ then the solution $b_{\a,\e,\d,n}$ of (\ref{ODE-1})-(\ref{ODE-3}) with $\zeta^{\e,\d}_n$ in place of $\zeta^{\e,\d}$ is an approximation of $b_{\a,\e,\d}$. Again, as 
$$\|b'_{\a,\e,\d,n}- b'_{\a,\e,\d}\|_{H^1(I)} \leq \|b_{\a,\e,\d,n}-b_{\a,\e,\d}\|_{H^2(I)},$$
executing Step (iii) $b'_{\a,\e,\d,n}$ is our desired approximate regularized solution.
Hence, if $j^{\d}$ and $g^{\e}$ are such that either $\zeta^{\e,\d}$ or $\zeta^{\e,\d}_n$ is in $R(T^{\a}_2)\cap C(I)$, then we have a stable approximate solution. Thus in this case we obtain a stable approximate solution to Problem (P) using steps among which the most critical one turns out to be that of solving an ODE.

\section{Appendix}\label{sec-7}
\bl
Let $J$ be a closed interval in $\R$. Then,
\beq \label{ap-1} \|y\|_{L^{\infty}(J)} \leq C_J \|y\|_{H^1(J)},\eeq
where $$C_J = C \max\{ 3, (2|J|+1)\}.$$
In particular, for any interval $J'$ contained in $J$,
\beq \label{ap-2} \|y\|_{L^{\infty}(J')} \leq C_J \|y\|_{H^1(J')}.\eeq
\el

\bpf
Let $J=[c,d]$ for some $c<d$. Let $\tilde{c}, \tilde{d} \in \R$ be such that $\tilde{c} < c$, $d< \tilde{d}$ and
\beq \label{ap-1-1} \max \{(c-\tilde{c}), (\tilde{d}-d) \} < (d-c).\eeq
Then, let us define the function
$$\tilde{y}(t)=   \left\{
\begin{array}{ll}
       0 & t\in \R\setminus[\tilde{c}, \tilde{d}] \\
       y(c)\left(\frac{t-\tilde{c}}{c-\tilde{c}}\right)  & t\in [\tilde{c},c]  \\
       y(t) & t \in J \\
       y(d)\left(\frac{\tilde{d}-t}{\tilde{d}-d}\right)  & t\in [d,\tilde{d}]

\end{array}
\right. $$
Then, it can be seen that $\tilde{y} \in H^1(\R)$ and
\beq \label{ap-1-2} \|\tilde{y}\|^2_{L^2([\tilde{c}, \tilde{d}])} = \frac{(y(c))^2}{(c-\tilde{c})^2}\int_{\tilde{c}}^c (t-\tilde{c})^2 dt + \|y\|^2_{L^2([c,d])} + \frac{(y(d))^2}{(\tilde{d}-d)^2}\int_d^{\tilde{d}} (\tilde{d}-t)^2 dt. \eeq
Now,
\beq \label{ap-1-3} \frac{(y(c))^2}{(c-\tilde{c})^2}\int_{\tilde{c}}^c (t-\tilde{c})^2 dt = \frac{(y(c))^2}{3(c-\tilde{c})^2} (c-\tilde{c})^3= \frac{(y(c))^2}{3} (c-\tilde{c})\eeq
and
\beq \label{ap-1-4} \frac{(y(d))^2}{(\tilde{d}-d)^2}\int_d^{\tilde{d}} (\tilde{d}-t)^2 dt = \frac{(y(d))^2}{3(\tilde{d}-d)^2}(\tilde{d}-d)^3= \frac{(y(d))^2}{3}(\tilde{d}-d).\eeq
By the fundamental theorem of calculus, for any $t\in [c,d]$,
$y(c)= -\int_c^t y'(s) ds + y(t),$
which implies,
$$|y(c)|^2= |y(t) -\int_c^t y'(s) ds |^2 \leq 2(|y(t)|^2 + |\int_c^t y'(s) ds |^2 ).$$
Hence, using Schwartz inequality as we have
$$|\int_c^t y'(s) ds |^2  \leq \left( \int_c^t |y'(s)|ds \right)^2 \leq (t-c)\|y'\|^2_{L^2([c,d])},$$
$$|y(c)|^2 \leq 2(|y(t)|^2 + (t-c)\|y'\|^2_{L^2([c,d])} )$$ holds.
This implies
$$|y(c)|^2(d-c)=\int_c^d |y(c)|^2 dt \leq 2\left( \int_c^d |y(t)|^2 dt + \|y'\|^2_{L^2([c,d])} \int_c^d (t-c) dt \right).$$
Thus,
\beq \label{ap-1-5}|y(c)|^2(d-c) \leq 2(\|y\|^2_{L^2([c,d])} + (d-c)^2 \|y'\|^2_{L^2([c,d])}).\eeq
Again, by the fundamental theorem of calculus, for any $t\in [c,d]$,
$y(d)= \int_t^d y'(s) ds + y(t),$
which implies,
$$|y(d)|^2= |y(t) +\int_t^d y'(s) ds |^2 \leq 2(|y(t)|^2 + |\int_t^d y'(s) ds |^2 ).$$
Hence, using Schwartz inequality as we have
$$|\int_t^d y'(s) ds |^2  \leq \left( \int_t^d |y'(s)|ds \right)^2 \leq (d-t)\|y'\|^2_{L^2([c,d])},$$
$$|y(d)|^2 \leq 2(|y(t)|^2 + (d-t)\|y'\|^2_{L^2([c,d])} )$$ holds.
This implies
$$|y(d)|^2(d-c)=\int_c^d |y(d)|^2 dt \leq 2\left( \int_c^d |y(t)|^2 dt + \|y'\|^2_{L^2([c,d])} \int_c^d (d-t) dt \right).$$
Thus,
\beq \label{ap-1-6}|y(d)|^2(d-c) \leq 2(\|y\|^2_{L^2([c,d])} + (d-c)^2 \|y'\|^2_{L^2([c,d])}).\eeq
Hence, combining (\ref{ap-1-2}), (\ref{ap-1-3}), (\ref{ap-1-4}), (\ref{ap-1-5}) and (\ref{ap-1-6}), we obtain
\begin{eqnarray}
\label{ap-1-7} \|\tilde{y}\|^2_{L^2([\tilde{c}, \tilde{d}])} 
&\leq & \frac{4}{3}(\|y\|^2_{L^2([c,d])} + (d-c)^2 \|y'\|^2_{L^2([c,d])}) + \|y\|^2_{L^2([c,d])} \\
\nonumber&\leq& \frac{7}{3} \|y\|_{L^2([c,d])} + \frac{4}{3}(d-c)^2 \|y'\|^2_{L^2([c,d])}.
\end{eqnarray}
Now,
$$\tilde{y}'(t)=   \left\{
\begin{array}{ll}
       0 & t\in \R\setminus[\tilde{c}, \tilde{d}] \\
       y(c) & t\in [\tilde{c},c]  \\
       y'(t) & t \in J \\
      -y(d) & t\in [d,\tilde{d}]

\end{array}
\right. $$
Hence,
\begin{eqnarray}
\label{ap-1-8} \|\tilde{y}'\|^2_{L^2([\tilde{c}, \tilde{d}])} 
&\leq & \int_{\tilde{c}}^c (y(c))^2 dt + \|y'\|^2_{L^2([c,d])} + \int_d^{\tilde{d}}(y(d))^2 dt\\
\nonumber &\leq & (y(c))^2(c-\tilde{c}) + \|y'\|^2_{L^2([c,d])} + (y(d))^2(\tilde{d}-d).
\end{eqnarray}
Thus, from (\ref{ap-1-5}) and (\ref{ap-1-6}), we obtain
\begin{eqnarray}
\label{ap-1-9} \|\tilde{y}'\|^2_{L^2([\tilde{c}, \tilde{d}])} 
&\leq & 2(\|y\|^2_{L^2([c,d])} + (d-c)^2 \|y'\|^2_{L^2([c,d])}) + \|y\|^2_{L^2([c,d])} + \|y'\|^2_{L^2([c,d])}\\
\nonumber &\leq& 3 \|y\|^2_{L^2([c,d])} + (2(d-c)^2 +1)\|y'\|^2_{L^2([c,d])}
\end{eqnarray}
and  from (\ref{ap-1-7}) and (\ref{ap-1-9}), we obtain
\beqarray
\|\tilde{y}\|_{H^1([\tilde{c}, \tilde{d}])} 
&\leq & 3\|y\|_{L^2([c,d])} + \left(\sqrt{\frac{10}{3}(d-c)^2+1}\right)\|y'\|_{L^2([c,d])} \\
&\leq& 3\|y\|_{L^2([c,d])} + (2(d-c)+1)\|y'\|_{L^2([c,d])}.\eeqarray
Hence,
\beq \label{ap-1-10} \|\tilde{y}\|_{H^1([\tilde{c}, \tilde{d}])} \leq \max\{ 3, (2(d-c)+1)\} \|y\|_{H^1([c,d])}.\eeq
Since, $H^1(\R)$ is continuously imbedded in $C(\R)\cap L^\infty(\R)$ (cf. \cite{Kes}), there exists $C>0$ such that
$\|\tilde{y}\|_{L^{\infty}(\R)} \leq C \|\tilde{y}\|_{H^1(\R)}, $
so that 
$$\|y\|_{L^{\infty}([c,d])} \leq \|\tilde{y}\|_{L^{\infty}(\R)} \leq C \|\tilde{y}\|_{H^1(\R)} = C \|\tilde{y}\|_{H^1([\tilde{c}, \tilde{d}])}.$$
Hence, by (\ref{ap-1-10})
$\|y\|_{L^{\infty}([c,d])} \leq C \max\{ 3, (2(d-c)+1)\} \|y\|_{H^1([c,d])}.$
\epf

\end{document}